\documentclass[12pt]{amsart}
\usepackage{latexsym}
\usepackage{amssymb}
\usepackage{amsfonts}
\usepackage{epsfig}
\usepackage{psfrag}
\usepackage{amsmath}
\usepackage{stackrel}
\usepackage{amscd}
\usepackage{psfrag}
\usepackage[all]{xy}

\newtheorem{defn0}{Definition}[section]
\newtheorem{prop0}[defn0]{Proposition}
\newtheorem{thm0}[defn0]{Theorem}
\newtheorem{lemma0}[defn0]{Lemma}
\newtheorem{corollary0}[defn0]{Corollary}
\newtheorem{example0}[defn0]{Example}
\newtheorem{conjecture0}[defn0]{Conjecture}
\newtheorem{notation0}[defn0]{Notation}
\newtheorem{remark0}[defn0]{Remark}
\newtheorem{quest0}[defn0]{Question}

\newenvironment{defn}{\begin{defn0}\rm}{\end{defn0}}
\newenvironment{prop}{\begin{prop0}}{\end{prop0}}
\newenvironment{thm}{\begin{thm0}}{\end{thm0}}
\newenvironment{lem}{\begin{lemma0}}{\end{lemma0}}
\newenvironment{cor}{\begin{corollary0}}{\end{corollary0}}
\newenvironment{example}{\begin{example0}\rm}{\end{example0}}
\newenvironment{rem}{\begin{remark0}\rm}{\end{remark0}}
\newenvironment{quest}{\begin{quest0}\rm}{\end{quest0}}

\newenvironment{notation}{\begin{notation0}\rm}{\end{notation0}}

\renewenvironment{proof}{\noindent {\textsc{Proof.}}}{$\square$ \vspace{3mm}}

\newenvironment{proofof}[1]{\noindent {\textsc{Proof of {#1}.}}}{$\square$ \vspace{3mm}}

\newcommand{\im}{\mathrm{Im}}

\newcommand{\T}{{\mathcal T}}

\newcommand{\m}{{\mathbf m}}
\renewcommand{\o}{{\mathbf o}}

\renewcommand{\sf}[2]{Tf_{#1}(#2)}

\newcommand{\x}{\texttt{X}}
\renewcommand{\a}{\texttt{A}}
\renewcommand{\c}{\texttt{C}}
\newcommand{\g}{\texttt{G}}
\renewcommand{\t}{\texttt{T}}
\renewcommand{\u}{\texttt{u}}
\renewcommand{\v}{\texttt{v}}

\newcommand{\I}{\mathcal{I}}
\newcommand{\id}{\textbf{1}}

\newcommand{\CC}{\mathbb{C}}

\newcommand{\NN}{\mathbb{N}}
\newcommand{\LL}{\mathcal{L}}

\newcommand{\lra}{\longrightarrow}
\newcommand{\ra}{\rightarrow}
\renewcommand{\ker}{\operatorname{Ker}\,}

\newcommand{\srel}{\stackrel}

\newif\ifprivate
\privatetrue

\def\???{\ifprivate {\bf {???}} \marginpar{{\Huge {\bf ?}}}
\else \fi}

 \DeclareMathOperator{\Hom}{Hom}

\DeclareMathOperator{\brk}{\mathbf{rk}\, }
\DeclareMathOperator{\rk}{rk\, } \DeclareMathOperator{\rep}{Par}
\DeclareMathOperator{\Int}{Int}

\numberwithin{equation}{section}

\begin{document}

\title[Relevant phylogenetic invariants...]{Relevant phylogenetic invariants of evolutionary models}

\author{Marta Casanellas}
\author{Jes\'{u}s Fern\'{a}ndez-S\'{a}nchez} \thanks{Both authors are partially supported by Ministerio de Educaci\'on y Ciencia,  MTM2009-14163-C02-02, and Generalitat de Catalunya, 2009 SGR 1284.}


\maketitle
\begin{abstract} Recently there have been several attempts to provide a
whole set of generators of the ideal of the algebraic variety
associated to a phylogenetic tree evolving under an algebraic
model. These algebraic varieties have been proven to be useful in
phylogenetics. In this paper we prove that, for phylogenetic
reconstruction purposes, it is enough to consider generators
coming from the edges of the tree, the so-called edge invariants.
This is the algebraic analogous to Buneman's Splits Equivalence
Theorem. The interest of this result relies on its potential
applications in phylogenetics for the widely used evolutionary
models such as Jukes-Cantor, Kimura 2 and 3 parameters, and
General Markov models.
\end{abstract}

\section{Introduction}
Algebraic evolutionary models and the algebraic varieties
associated to a tree evolving under these models have been an
interdisciplinary area of research with successful results in the
last five years.  The use of polynomials in phylogenetic
reconstruction was first introduced by biologists Cavender and
Felsenstein \cite{Cavender87} and Lake \cite{Lake1987}. Because of
their interest in phylogenetics, there have been several attempts
to provide a set of generators of the ideal of these algebraic
varieties (see for example \cite{Allman2004b},
\cite{Sturmfels2005}, \cite{Draisma}, \cite{CS}).
On the other hand, the authors of this paper have proven in
\cite{CFS} that these generators can be successfully used in
phylogenetic reconstruction. In other words, methods based in
algebraic geometry can lead to the inference of the phylogenetic
tree of current biological species. As we already did in
\cite{CFS2}, our aim in the present paper  is  to address  again
the study of
these algebraic varieties towards 
their real applications in phylogenetics. 

Algebraic evolutionary models include the algebraic version of
widely used models in biology such as Jukes-Cantor model
\cite{JC69}, Kimura 2 and 3 parameters model (cf.
\cite{Kimura1980}, \cite{Kimura1981}) and the general Markov model
(cf. \cite{barryhartigina87}). These models belong to what Draisma
and Kuttler call \textit{equivariant models} in \cite{Draisma} (see section 2 for
the precise definition). Following ideas of Allman and Rhodes and
using representation theory, Draisma and Kuttler have recently
given an algorithm to obtain the generators of the ideal of the
algebraic varieties associated to a tree of $n$ species evolving
under an equivariant model from the generators of the ideal
associated to a tree of 3 species and certain minors of matrices
(the so-called \textit{edge invariants}). Nevertheless, a set of
generators for trees of 3 species is not known for certain models
such as the general Markov model (this is the so-called Salmon
Conjecture) or the strand symmetric model (see \cite{CS}).
Therefore, a complete list of generators for a tree of $n$ species
evolving under these models cannot be given at this point.

The goal of this paper is to prove that, whereas mathematically
speaking it is interesting to know a set of generators of the ideal
of these varieties, for biological purposes it is enough to
consider certain generators. More precisely, the edge
invariants mentioned above suffice to reconstruct the phylogenetic
tree of any number of species (see the Theorem in the next page or Theorem \ref{teophyloinvar}). This is a
natural result if one thinks of the combinatorics result of
Buneman that says that a tree can be recovered if one knows the
set of splits on the set of leaves induced by its edges (cf.
\cite{Buneman}, \cite[Theorem 2.35]{ASCB2005}, see also Theorem
\ref{splits} below).

Our inspiration goes back to the work \cite{Felsenstein1991} of
biologist Joe Felsenstein who calls  \emph{phylogenetic
invariants} to those polynomial expressions that vanish on the
expected frequencies of any sequences arising from one tree
topology but are non zero for at least one tree of another
topology. A tree topology in this setting is the topology of the
tree graph labelled at the leaves with the name of the species.
Algebraically speaking, he calls \emph{phylogenetic invariants} to
those elements of the ideal associated to a phylogenetic tree that
allow to distinguish it from other tree topologies. In the
mathematical context, the name phylogenetic invariants has usually
been given to all elements of the ideal, see for instance the work
of Allman and Rhodes \cite{Allman2004b}. We want to go back to the
original meaning of phylogenetic invariants because our focus is
devoted to the applications of algebraic geometry in the
reconstruction of the tree topology of current species. Therefore,
we are mainly interested in precisely those elements of the ideal
that provide information for phylogenetic reconstruction purposes;
in other words, we are interested in \textit{phylogenetic}
invariants (i.e polynomials in the ideal of one tree topology of
$n$ species but not in the ideal of all other tree topologies on
the same number of species) and the word \textit{invariants} alone
shall mean any element of the ideal. In colloquial language the
main result of this paper is that, for phylogenetic reconstruction
purposes, the relevant phylogenetic invariants are the edge
invariants mentioned above.

As our aim is to study these varieties regarding their
applications in biology, let us roughly explain here how does
algebraic geometry interfere with phylogenetic reconstruction. Let
$n$ be a number of biological species and assume that we are given
an alignment of DNA sequences corresponding to them (the
definition of alignment is rather technical but it refers to a
collection of $n$-tuples in $\{\texttt{A,C,G,T}\}^n$ that will be
also called columns of the alignment). Each column stands for sites in the $n$ DNA sequences that 
have evolved from the same nucleotide in the common ancestor. We assume that these
species are leaves of a phylogenetic tree $T$ evolving under a
probabilistic model $\mathcal{M}$ (in this paper we will only
consider equivariant models, see Definition 2.4 for the precise
definition).
It is usual to assume as well that all columns of the alignment
behave independently and identically (i.e. all sites of the DNA
sequences of these species evolve in the same way and
independently of the other sites). Associated to this model
$\mathcal{M}$ there is a parameterization map $\Psi_T$ giving the
joint distribution of states $\texttt{A,C,G,T}$ at the leaves of
$T$  as polynomial functions of continuous parameters. Therefore,
as an alignment of DNA sequences evolving under this model on a
tree $T$ is a collection of observations of states at the leaves,
it corresponds to a point in the image of this parameterization
map. The algebraic variety $V_{\mathcal{M}}(T)$ associated to $T$
is the closure of this image (see Definition 2.7). In the real
life, alignments are not points of $V_{\mathcal{M}}(T)$ but they
are \textit{close} to $V_{\mathcal{M}}(T)$ if the model reasonably
fits the data. Therefore the idea behind phylogenetic algebraic
geometry is to use the ideal of $V_{\mathcal{M}}(T)$ in order to
infer the tree topology $T$. See \cite{CGS} for an algorithm of phylogenetic reconstruction based on the generators of this ideal and \cite{CFS} for tests of it on simulated data. 

Up to now, all attempts have focused on giving a whole set of
generators of $I(V_{\mathcal{M}}(T))$ but our approach is more
practical. As biologists assume that the model $\mathcal{M}$ fits
the data, the point given by an alignment is therefore assumed to
be close to the union of all varieties $V_{\mathcal{M}}(T)$ for
trees of $n$ species evolving under model $\mathcal{M}$.
Henceforth, we only need to know how is a particular variety
$V_{\mathcal{M}}(T_0)$ defined inside $\cup_T V_{\mathcal{M}}(T)$
where the union runs over all trivalent tree topologies $T$ of $n$
species. In this algebraic geometry context our main result
(Theorem \ref{teophyloinvar}) can be summarized in the following
way.

\vspace{2mm}

\textbf{Theorem.} \textit{Let $\mathcal{T}$ be the set of
trivalent tree topologies on $n$ leaves and let $\mathcal{M}$ be
an equivariant model. For each tree topology $T \in \mathcal{T}$
there exists an open set $U_T$ such that if $p$ belongs to
$\cup_{T \in \mathcal{T}} U_T$, then $p$ belongs to a particular
variety $V_{\mathcal{M}}(T_0)$ if and only if $p$ belongs to the
zero set of edge invariants of $T_0$.}

\vspace{2mm}

This result has also other consequences in phylogenetics. For
instance, it says that edge invariants should not be used for
model fitting tests  (see \cite{Garcia} for an algebraic
introduction to the subject) or for the study of identifiability
of continuous parameters (see \cite{Allman2008} for an explanation
of these terminology) of the model because they are indeed
phylogenetic invariants. Instead, they should be used in
discussing the identifiability of tree topology of such models
(see Corollary \ref{identifiability}) as it was already done by
Allman and Rhodes in \cite{Allman2006}. We also find invariants
(not phylogenetic invariants) that could potentially be used for
model fitting tests, that is, linear polynomials that can be used
for choosing the evolutionary model that best fits the data (see
Remark \ref{modelfitting}).

Moreover, our main theorem allows one to give the exact degrees of
those generators relevant in phylogenetics (see Corollary
\ref{degrees}), whereas the degrees of a whole set of generators
for the general Markov or strand symmetric models are still
unknown. It is worth highlighting that these degrees can be
computed by just knowing the model we are interested in, and they
do not depend on the topology or the number of leaves we are
considering.

Here we outline the structure of the paper. In section 2 we adapt
the setting and notation of \cite{Draisma} to our convenience. As
well, we prove and recall basic facts of group representation
theory for those non-familiarized readers. Section 3 is devoted to
prove a technical result that will be the key in the proof of our
main theorem. Roughly speaking this result proves that edge
invariants are indeed \textit{phylogenetic} invariants for any
equivariant model. This was already known  for the general Markov
model by Allman, Rhodes (see for instance \cite{Allman2006}) and
Eriksson \cite{Eriksson05} but it is new for the remaining
equivariant models. The proof relies on providing a formula for
the rank of the flattening of the tensor $\Psi_T$ along any
bipartition of the set of leaves. In section 4 we prove Theorem
\ref{teophyloinvar}, our main result. In the last section we
provide an exhaustive collection of examples on how to compute the
required edge invariants for the most used evolutionary models:
Jukes-Cantor, Kimura 2 and 3 parameters, strand symmetric and
general Markov model. We compute them explicitly for quartet
trees. It is our aim to make this section clear enough for
biomathematicians so that, for example, we relate invariants used
by biologist like Lake (see \cite{Lake1987}) to the more technical
definition of edge invariants (see the end of subsection 5.5). We
also connect our edge invariants to Fourier coordinates that are
more familiar to those readers used to group-based models. In
particular, the reader can visualize what are the Fourier
coordinates that are actually interesting in biology as not all of
them are needed for phylogenetic reconstruction. This section is
also a useful illustration of technical definitions given in
sections 2 and 3 so it is a good idea to combine the reading of
both sections with section 5.

\vspace{2mm}

\textbf{Acknowledgments:} We would like to thank Josep Elgueta and Jeremy Sumner for useful comments on group
representation theory.
%

\section{Preliminaries}
A \emph{tree} is a connected finite graph without cycles,
consisting of vertices and edges. Given a tree $T$, we write
$V(T)$ and $E(T)$ for the set of vertices and edges of $T$.  The
\emph{degree} of a vertex is the number of edges incident on it.
The set $V(T)$ splits into the set of \emph{leaves} $L(T)$
(vertices of degree one)  and the set of interior vertices
$\Int(T)$: $V(T)=L(T)\cup \Int(T)$.
One says that a tree is \emph{trivalent} if each vertex in
$\Int(T)$ has degree 3.  A \emph{tree topology} is the topological
class of a tree where every leaf has been labelled.
Given a subset $L$ of $L(T)$, the \emph{subtree induced by} $L$ is
just the smallest tree composed of the edges and vertices of $T$
in any path connecting two leaves in $L$.


%

\vspace{2mm} Given an ordered set $B=\{b_1,b_2,\ldots,b_k\}$, we
define $W=\langle B \rangle_{\CC}$ as the $\CC$-vector space
generated by the elements of $B$.
For biological applications,  the most common values of $k$ are 2,
4 or 20 (for example,
$B=\{\texttt{A},\texttt{C},\texttt{G},\texttt{T}\})$.
Now, given a subgroup $G$ of the group $\mathfrak{S}_k$ of
permutations of $k$ elements, we consider the restriction to $G$ of
the natural linear representation
\begin{eqnarray*}
  \rho:\mathfrak{S}_k\ra GL(W)
\end{eqnarray*}
given by the permutation of the elements of $B$. %
This representation induces a $G$-module structure on $W$ by
taking
\begin{eqnarray*}
g\cdot u:=\rho(g)(u)\in W.
\end{eqnarray*}
In fact, $\rho$ induces a $G$-module structure on any tensor power
of $W$, say  $\otimes^l W:=W\otimes \ldots \otimes W$, by taking
\begin{eqnarray}\label{mult_action}
g\cdot \left (u_{1}\otimes \ldots \otimes u_{l}\right ):=g\cdot
u_{1} \otimes \ldots \otimes g\cdot u_{l}.
\end{eqnarray}
Henceforth, any tensor power of $W$ will be implicitly considered
as a $G$-module with this action.

\vspace{2mm} From now on, we fix an ordered set
$B=\{b_1,b_2,\ldots,b_k\}$, $W=\langle B \rangle_{\mathbb{C}}$ and
a subgroup $G\subset \mathfrak{S}_k$ acting on $W$ as above.

\begin{defn}
A \emph{phylogenetic tree on} $(G,W)$ is a tree where every vertex
$p$ has a $\CC$-vector space $W_p\cong W$ associated to it,
regarded as a representation of $G$ via the map $\rho$ defined
above.
\end{defn}

%
%
%

%

\noindent \textbf{Notation. }The scalar product with orthonormal basis $B_p$
will be denoted by $(. \mid  .)_p$.  This gives a canonical
isomorphism from $W_p$ to $W_p^{\ast}$.

\vspace{2mm}
Notice that the scalar product $(. \mid  .)_p$  is $G$-invariant,
that is, $(g\cdot u\mid g \cdot v)_p=(u\mid v)_p$ for every
$u,v\in W_p$ and any $g\in G$.

\begin{defn}
Given a phylogenetic tree $T$ on $(G,W)$, a $T$-\emph{tensor} is any element of
\[\LL(T):=\otimes_{p\in L(T)} W_p.\]
A \emph{$G$-tensor} on $T$ is a $T$-tensor  invariant by the action defined in (\ref{mult_action}). The set of $G$-tensors will be denoted by  $\LL(T)^G$.
\end{defn}

From now on, if $l>0$ we write $\otimes^l W=W\otimes  \srel{l}
\ldots \otimes  W$.
We denote by $B(\otimes^l W)$ the basis of $\otimes^l W$ given by \[\{u_{i_1}\otimes  \ldots \otimes u_{i_l}\mid u_{i_j}\in B\}.\]
This is an othonormal basis with respect to the scalar product of $\otimes^l W$ given by $(\otimes_p u_p\mid \otimes_p v_p)=\prod_p (u_p\mid v_p)$. If $L\subset L(T)$ is a subset of $L(T)$ and $l=\sharp L$, then we shall use the notation $\otimes_L W$ for the space $\otimes_{p\in L} W_p\cong \otimes^l W$.

\begin{defn}Let $T$ be a phylogenetic tree on $(G,W)$ and assume that a distinguished vertex of $T$ (the \emph{root}) is given, inducing an orientation in all the edges of~$T$: write $e_0$ and $e_1$ for the origin and final vertices of the edge  $e$, respectively.
A  $G$-\emph{evolutionary presentation}\footnote{Notice that
\emph{evolutionary presentations} are called
\emph{representations} in \cite{Draisma}. We prefer this
terminology to avoid confusion with representation theory.} of $T$
is a collection of tensors  $\{A_{e_0,e_1}\}_{e\in E(T)}$ where each
$A_{e_0,e_1}$ is a $G$-invariant element of the $G$-module $W_{e_0} \otimes W_{e_1}$.
The space of $G$-invariant elements of $W_{e_0} \otimes W_{e_1}$ is denoted by $\left ( W_{e_0} \otimes W_{e_1}\right )^G$.

%
If another root (orientation) on $T$ is considered, inducing the
opposite orientation on some edge $e\in E(T)$, we define
$A_{e_1,e_0}:=A_{e_0,e_1}^t$, where $.^t$ is the natural
isomorphism $\left ( W_{e_0} \otimes W_{e_1} \right )^G\cong \left ( W_{e_1} \otimes
W_{e_0} \right )^G$. We will often identify $\Hom_G(W_{e_0},W_{e_1})$ with $\left ( W_{e_0}\otimes W_{e_1}\right )^G$ via $W_{e_0}^{\ast}\cong W_{e_0}$.
With this convention, $G$-evolutionary presentations on a tree do not
depend on the orientation chosen.
The space of all $G$-evolutionary presentations of $T$ is the
parameter space denoted by $\rep_G(T)=\prod_{e\in E(T)} \left ( W_{e_0}\otimes W_{e_1}\right )^G$.
Notice that a $G$-evolutionary presentation of $T$ induces  by
restriction a $G$-evolutionary presentation of any subtree of $T$.

The space  $\rep_G(T)$, as well as $\LL(T)$ and
$\LL(T)^G$, are irreducible affine spaces with their Zariski
topology.
\end{defn}


\begin{defn}
An \emph{equivariant model} of evolution is a pair $(G,W)$ as
above, $W=\langle b_1, \dots,b_k \rangle$, $G \subset
\mathfrak{S}_k$. Trees evolving under this equivariant model  are
phylogenetic trees on $(G,W)$ together with the space of
$G$-evolutionary presentations.

Equivariant models of evolution include the general Markov model \cite{barryhartigina87} when $G=\{\mathrm{id}\}$, the strand symmetric model \cite{CS} when $G=\langle (\a\t)(\c\g)\rangle$, and the algebraic versions of Kimura 3-parameters \cite{Kimura1981} ($G=\langle (\a\c)(\g\t),(\a\g)(\c\t)\rangle$), Kimura 2-parameters \cite{Kimura1980} ($G=\langle (\a\c\g\t),(\a\g)\rangle$) and Jukes-Cantor models \cite{JC69} ($G=\mathfrak{S}_4$). We derive the reader to section 5 for specific computations with these models.

\end{defn}

%

Following \cite{Allman2004b} and \cite{Draisma} we present now a
fundamental operation $*$ on phylogenetic trees, $G$-evolutionary
presentations and $T$-tensors. To this aim, we first introduce a bilinear operation $\langle \cdot \mid \cdot \rangle$ between tensors induced by the bilinear form $(\cdot \mid \cdot)$ on $W$.
Let $X$ and $Y$ be two finite sets of indices with $Z=X\cap Y\neq \emptyset$, and such that every $p$ in $X$ or $Y$ has associated a vector space $W_p\cong W$ to it.
Define
\begin{eqnarray*}\label{forflat}
\langle .\mid  .\rangle:  \otimes_{X}W \times \otimes_{Y}W  &
\rightarrow & \otimes_{X\cup Y \setminus Z}
W \nonumber \\
(\otimes_{p\in X} v_p,\otimes_{p\in Y} u_p) & \mapsto & (\otimes_{p\in Z}v_p\mid \otimes_{p\in Z}u_p) \left ( (\otimes_{p\in X\setminus Z} v_p) \otimes
(\otimes_{p\in Y\setminus Z} u_p) \right )
\end{eqnarray*}

\vspace{2mm}

Now, we define the $\ast$ operation:

\vspace{2mm}
\noindent \textbf{$*$ for trees:} Given $l$ spaced trees $T_1, . .
. , T_l$ whose vertex sets only share a common leaf $q$ with
common space $W_q$ and common basis $B_q$,
we construct a new spaced tree $*_iT_i$ obtained by gluing the
$T_i$'s' along $q$; the space at a vertex of $*_iT_i$ coming from
$T_j$ is just the space attached to it in $T_j$, with the same
distinguished basis.

\noindent \textbf{$*$ for $G$-evolutionary presentations:} Given
$G$-evolutionary presentations $A_i\in  \rep(T_i)$ for $i = 1, . . . ,
l$, we denote by $*_iA_i$ the $G$-evolutionary presentation of
$*_iT_i$ built up from the $A_i$.

\noindent \textbf{$*$ for tensors:} Now let $\psi_i$ be a
$T_i$-tensor, for all $i$. Then we obtain a $T$-tensor as follows:
\[*_i\psi_i := \sum_{b\in B_q}
\otimes_i \langle b \mid \psi_i \rangle.\]
Although this $*$ operator is not a binary operator extended to
several factors, when convenient we will write $T_1 * \ldots *
T_l$ for $*_iT_i$ and $\psi_1 *\ldots  * \psi_l$ for $*_i \psi_i$.

\vspace{3mm}

\begin{notation}\label{new}
A slightly more general $\ast$-operation will be needed in forthcoming section~3.
Given $\varphi_A\in (\otimes_{X} W)^G$ and $\varphi_2\in (\otimes_{Y} W)^G$, define
\begin{eqnarray*}
\varphi_1 \stackrel[Z]{}{\ast} \varphi_2=\sum_{b\in B(\otimes_{Z} W)}\langle
\varphi_1\mid b\rangle \otimes \langle \varphi_2 \mid b\rangle \in
(\otimes_{X\cup Y \setminus Z} W)^G.
\end{eqnarray*}
Clearly, if $T_1$ and $T_2$ are two phylogenetic trees that share a common leaf $q$, then this definition agrees with the $\ast$-operation defined above.
\end{notation}


\vspace{3mm}

Now we describe a basic procedure that allows us to associate a
$T$-tensor to any $G$-evolutionary presentation of~$T$. We proceed
inductively on the number of edges to define $\Psi_T : \rep(T ) \rightarrow \LL(T )$. Let
$A \in \rep(T )$. First, if $T$ has a single edge ${p, q}$, then
$\Psi_T(A) := A_{qp}$, is an element of $\LL(T) = W_q \otimes
W_p$. If $T$ has more than one edge, then let $q$ be any internal
vertex of $T$.
Two vertices $p,q\in T$ are
\emph{adjacent} if they are joined by an edge; in this case, we write
$p\sim q$.
We can then write $T = *_{p\sim q}T_p$, where $T_p$
is the branch of $T$ around $q$ containing $p$, constructed by
taking the connected component of $T\setminus \{q\}$ containing
$p$, and reattaching $q$ to $p$. The $G$-evolutionary presentation $A$
induces $G$-evolutionary presentations $A_p$ of the $T_p$, and by
induction $\Psi_{T_p}(A_p)$ has been defined. We now set
\[\Psi_T (A) := *_{p\sim q}\Psi_{T_p} (A_p).\]
This definition is independent of the choice of $q$ and the
formula is also valid if $q$ is actually a leaf (see
\cite{Draisma} for details).  Moreover, we have that the map
${\Psi_T:\rep_G(T ) \rightarrow \LL(T )}$ is $G$-equivariant (see
\cite[Lemma 5.1]{Draisma}).

\begin{rem}\label{RmkZariski}
Notice that the above map $\Psi_T : \rep_G(T) \rightarrow \LL(T)^G$
is a continuous map in the Zariski topology.
\end{rem}

\begin{defn}
The algebraic variety associated to a phylogenetic tree $T$ on $(G,W)$ is
\[V_G(T) := \overline{\left \{\Psi_T (A)\mid A\in  \rep_G(T ) \right \}} \subset \LL(T )\]where the closure is taken in the Zariski topology.

Notice that we have $V_G(T) \subset \LL(T)^G$. From now on, we
will consider $\LL(T)^G$ as the ambient space of $V_G(T)$ and
$\I(T)$ will be the ideal of this variety in the corresponding
coordinate ring. When the group is understood from the context, we will use the notation $V(T)$.
\end{defn}

\begin{rem}\label{modelfitting} The inclusion $\LL(T)^G \subseteq \LL(T)$ is defined by a
set of linear polynomials that are also \textit{invariants} of any
phylogenetic tree $T$ on $(G,W)$ (see the Introduction for the
explanation of the word invariants). Although they are not
\textit{phylogenetic invariants} because they vanish on $V_G(T)$
for any tree $T$, they might be interesting for choosing the model
$(G,W)$ that best fits the data. This application of invariants to
model fitting will be studied in a forthcoming paper.
\end{rem}

\begin{example}
If we consider $B=\{\a,\c,\g,\t\}$ and $G=\{\textrm{id}\}\subset \mathfrak{S}_4$, we obtain the general Markov model. In this case, $\left ( W_{e_0}\otimes W_{e_1}\right )^G=\left ( W_{e_0}\otimes W_{e_1}\right )$ and no restrictive conditions are imposed on the parameters of the model. Thus, a $G$-evolutionary presentation can be identified, by taking the basis $B$ in $W$ with a collection of matrices $\{A_e\}_{e\in E(T)}$ and the parameters of the model are the entries of these matrices. When these entries are real non-negative values and their columns sum to 1, they can be understood as the probabilities of substitution among the 4 nucleotides:
\small
\begin{eqnarray*}
A_e=\left ( \begin{array}{cccc}
P(\a\mid \a,e) & P(\a\mid \c,e) & P(\a\mid \g,e) & P(\a\mid \t,e) \\
P(\c\mid \a,e) & P(\c\mid \c,e) & P(\c\mid \g,e) & P(\c\mid \t,e) \\
P(\g\mid \a,e) & P(\g\mid \c,e) & P(\g\mid \g,e) & P(\g\mid \t,e) \\
P(\t\mid \a,e) & P(\t\mid \c,e) & P(\t\mid \g,e) & P(\t\mid \t,e)
\end{array} \right ).
\end{eqnarray*}
\normalsize Here $P(\x \mid \texttt{Y},e)$ is the conditional
probability that nucleotide \texttt{Y} at the parent species $e_0$
is being substituted along edge $e$ by nucleotide \texttt{X} at
its child species $e_1$. In our terminology introduced above,
$P(\x \mid \texttt{Y},e)$ is the coordinate of $A_e \in W_{e_0} \otimes W_{e_1}\cong W\otimes W $ corresponding to $ \texttt{Y} \otimes \x$. Given a tree $T$,
the $G$-equivariant map $\Psi_T$ is the parameterization that
associates
to each parameter set the
vector of expected pattern frequencies $p=(p_{\x_1 \x_2 \ldots
\x_n})_{\x_i\in B}$ (that is, $p_{\x_1 \x_2 \ldots \x_n}$ is the
probability of observing $\x_1 \x_2 \ldots \x_n$ at the leaves of
$T$). For example, if $T$ is a 4-leaf tree as in figure \ref{4-tree},
\begin{figure}\label{4leaves}
\begin{center}
\psfrag{e1}{$e(1)$}\psfrag{e2}{$e(2)$}
\psfrag{e3}{$e(3)$}\psfrag{e4}{$e(4)$}
\psfrag{e}{$e$}
\includegraphics[scale=0.8]{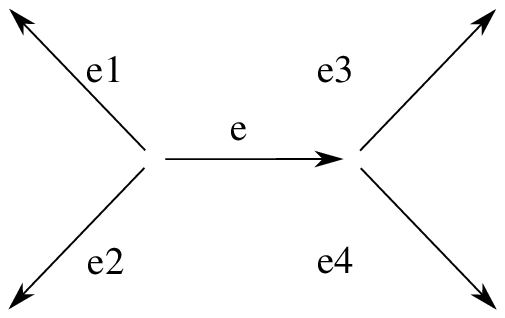}
\end{center}
\caption{\label{4-tree}}
\end{figure}
then
\begin{eqnarray*}
\Psi_T: \prod_{e\in E(T)} (W\otimes W)\cong \CC^{80} & \ra & \otimes^ 4 W \cong \CC^{256} \\
(A_e)_{e} & \mapsto & (p_{\a\a\a\a},p_{\a\a\a\c},\ldots,
p_{\t\t\t\t})
\end{eqnarray*}
and $p_{ \x_1 \x_2 \x_3 \x_4 }$ is the coordinate of $p \in
\LL(T)\cong\CC^{256}$ corresponding to the basis vector $\x_1
\otimes \x_2 \otimes \x_3 \otimes \x_4.$ In this case, the image
of $\Psi_T$ is given by
$$p_{ \x_1 \x_2 \x_3 \x_4 }=\sum_{\texttt{Y},\texttt{Z}}\pi_{\texttt{Y}}A_e(\texttt{Z},\texttt{Y})A_{e(1)}(\x_1,\texttt{Y})A_{e(2)}(\x_2,\texttt{Y})A_{e(3)}(\x_3,\texttt{Z})A_{e(4)}(\x_4,\texttt{Z}).$$
Here $\pi_{\texttt{Y}}$ is the probability of nucleotide
$\texttt{Y}$ occurring at the root node (see figure
\ref{4-tree}). Actually, in the original definition of $\Psi_T$
(see paragraph before Remark \ref{RmkZariski}) we gave a
reparameterization of $V_G(T)$ where we omit parameters
$\pi_{\texttt{Y}}$ for convenience.
\end{example}





\begin{defn}
Given a tree $T$, a \emph{bipartition} of the leaves of $T$ is a
decomposition $L(T)=L_1\cup L_2$ where $L_1\cup L_2=\emptyset $.
We denote it as $L_1\mid L_2$.
Notice that every edge $e$ of $T$ induces a bipartition $L_1 \mid
L_2$ of $L(T)$ by removing it; such a bipartition is called an
\emph{edge split} of $T$ and will be denoted by the same letter $e$.
\end{defn}

%

\subsection{Representation Theory}
\label{notation} We will make use of representation theory of
groups. A basic reference for this are the books \cite{Serre} and
\cite{FultHarr} and the reader is referred to them for definitions
and well-known facts.
%

From now on, write $\Omega_G=\{\omega_1, \dots, \omega_s\}$ for
the set of irreducible characters of $G$.
%
%
It is known that any two representations with the same character
are isomorphic (Corollary 2 of $\S$ 2 of \cite{Serre}). As a
consequence of this and  Schur's lemma (see $\S$2.2 of
\cite{Serre}) we obtain the following fundamental result in
representation theory:

\begin{lem}\label{alter_Schur}
Let $N_{\omega}, N_{\omega'}$ be the irreducible linear representations of $G$ with
associated  characters $\omega, \omega'\in \Omega_G$. If  $f:N_{\omega}
\rightarrow N_{\omega'}$ is a $G$-module homomorphism, and
\begin{itemize}
\item[(i)] if $\omega\neq \omega'$, then $f =0$;
\item[(ii)] if $\omega=\omega'$, then $f$ is a homothety.
\end{itemize}
In particular, $\Hom_G(N_{\omega}, N_{\omega})\cong \mathbb{C}$.
\end{lem}

For every irreducible character $\omega_t\in \Omega_G$, fix an
irreducible $G$-module $N_{\omega_t}$ with associated character
$\omega_t$. Then, for any $G$-module $V$, there exists a unique
decomposition of $V$ into isotypic components:
\begin{eqnarray}\label{dec}
V \cong  \oplus_{t=1}^s V[\omega_t]
\end{eqnarray}
where each $V[\omega_t]$ is isomorphic to $N_{\omega_t} \otimes
\CC^{m(\omega_t,V)}$ for some multiplicity $m(\omega_t,V)$,
$t=1,\ldots,s$.
We also have that if $V'$ is another representation of $G$, then
\begin{eqnarray}\label{schur}
\Hom_G(V, V') \cong
\oplus_{t=1}^{s}\Hom_{\CC}(\CC^{m(\omega_t,V)},
\CC^{m(\omega_t,V')})
\end{eqnarray}

Going back to our fixed vector space $W$, we already know that the
space $\otimes^l W$, $l>0$   is a $G$-representation as well and,
as such,
\[\otimes^l W \cong  \oplus_{t=1}^s N_{\omega_t} \otimes \CC^{m(\omega_t,\otimes^l W)}.\]
We will denote by $\m(l)$ the $s$-tuple
\[\m(l)=(m(\omega_1,\otimes^l W), \ldots,m(\omega_s,\otimes^l W)).\]
In particular, $\m(1)$ will be denoted by $\textbf{m}=(m_1,
\ldots,m_s)$.
Moreover, if $\chi$ denotes the associated character to the
representation $\rho:G\lra GL(W)$, the decomposition (\ref{dec})
above induces an equality of characters
\[\chi=\sum_{t=1}^s m_t \omega_t \quad m_t\in \mathbb{Z}.\]
If
$\mathbf{a}=(a_t)_{t=1,\ldots,s},\mathbf{b}=(b_t)_{t=1,\ldots,s}
\in \NN^s$, we write $\mathbf{a}\leq \mathbf{b}$ if $a_t\leq b_t$
for each  ${t=1, \dots,s}$. Similarly, $\min
\{\mathbf{a},\mathbf{b}\}$ is the $s$-tuple given by the minimum
of each entry.

\begin{lem}\label{max_m}
With this notation, we have $\m(l)\leq \m(l')$ if $l\leq l'$.
\end{lem}

\begin{proof}
We prove that $\m(l)\leq \m(l+1)$ for any $l$. First of all, we
show that if $\omega_1\in \Omega$ is the trivial character, then
$m_1\geq 1$. To this aim, notice that the vector $\sum_{b\in
B}b\in W$ is invariant by the action of any $g\in G$. In
particular, we have $\sum_{b\in B}b\in W[\omega_1]$ and so
$\omega_1$ does appear in the decomposition of $\chi$ with
non-zero coefficient.
Now, given $l>0$, write $\chi^l=\sum_t a_t\omega_t$. The claim
follows from the fact that the coefficient  of any irreducible
character of $G$, say $\omega_t$, in $\chi^{l+1}$ is just $m_1
a_t+\ldots \geq a_t$.
\end{proof}



\begin{notation}\label{ranknot}
If $A$ is a matrix in $M_{m,n}$, we say that $A$ has \emph{maximal
rank }if $\rk(A)=\min\{n,m\}$.
Following \cite{Allman2004b}, if $\m$, $\mathbf{n}$ are $s$-tuples
of positive integers we will use the notation $ M_{\m,\mathbf{n}}$
to denote the space $M_{m_1,n_1}\times \dots \times M_{m_s,n_s}$
and if $A=(A_1,\ldots,A_s)\in M_{\m,\mathbf{n}}$, we will write
\[\brk (A)=\left (\rk(A_1),\ldots,\rk(A_s)\right ).\]
Notice that $M_{\m,\mathbf{n}}$ can be understood as the subspace
of $M_{\sum m_t,\sum n_t}$ given by the block-diagonal matrices
with blocks of sizes $m_t\times n_t$. Then, $A\in
M_{\m,\mathbf{n}}$ has maximal rank  as a matrix in $M_{\sum
m_t,\sum n_t}$ if and only if  $\brk(A)=\min\{\m,\mathbf{n}\}$.
\end{notation}

\subsection{Flattenings and thin flattenings}

The following definitions will be crucial for our purposes.
\begin{defn}\label{SchurF}
Let $T$ be a phylogenetic tree on $(G,W)$ and let $L_1\mid L_2$ be a
bipartition of its leaves with  $l_1=\sharp L_1$ and $l_2=\sharp
L_2$. Let $\psi$ be a $G$-tensor on $T$.

The \emph{flattening of $\psi$ along} ${L_1\mid L_2}$, denoted by
$flat_{L_1\mid L_2}\psi$, is the image of $\psi$ via the
isomorphism
\[\LL(T)^G\cong \Hom_G\left  (\otimes_{L_1}W,\otimes_{L_2}W \right ).\]

The \emph{thin flattening of} $\psi$ \emph{along} ${L_1\mid L_2}$
is the  $s$-tuple of linear maps, denoted by $\sf{L_1\mid
L_2}{\psi}$, obtained from $flat_{L_1\mid L_2}\psi$ via the
isomorphism
\begin{eqnarray*}
\Hom_G(\otimes_{L_1}W,\otimes_{L_2}W) \cong  \bigoplus_{t=1}^s
\Hom_{\CC}(\CC^{\m(l_1)_t},\CC^{\m(l_2)_t}).
\end{eqnarray*}
\end{defn}

\begin{rem}
Notice that the composition of linear maps induce a composition of flattenings and thin flattenings.
Notice also that if $\psi\in \LL(T)^G$ and $L_1\mid L_2$ is a bipartition of $L(T)$, then
\[\left ( flat_{L_1\mid L_2} \psi \right )(u) = \langle \psi \mid u \rangle, \qquad  \forall u\in \otimes_{L_1} W \]where $\langle \cdot \mid \cdot \rangle$ is the operation defined in (\ref{forflat}).
\end{rem}

\begin{notation}\label{basis/flat}
If $\sf{L_1\mid L_2}{\psi}=(\psi_1,\psi_2,\ldots,\psi_s)$, we write
\[\brk \sf{L_1\mid L_2}{\psi}=(\rk
(\psi_1),\ldots,\rk (\psi_t)).\]We also denote
$ \rk \sf{L_1\mid L_2}{\psi} =\sum_{t=1}^s \rk (\psi_t) $
and call it the \emph{rank} of $\sf{L_1\mid L_2}{\psi}$.
Clearly, this definition is coherent with the usual definition of rank if we regard $\sf{L_1\mid L_2}{\psi}$ as a $\CC$-linear map $\CC^{\sum_t \m(l_1)_t}\ra \CC^{\sum_t \m(l_2)_t}$.
\end{notation}

\vspace{2mm}
The following easy lemma is left to the reader:

\begin{lem}\label{eq}
We have that $\rk flat_{L_1\mid L_2}\psi=\sum_{t=1}^s
\dim_{\mathbb{C}}N_{\omega_t} \rk (\psi_t)$. Moreover, the following are equivalent
\begin{itemize}
 \item[(i)] $\rk flat_{L_1\mid L_2}\psi$ is maximal;
\item[(ii)] $\rk (\psi_t)$   is maximal  $\forall t\in \{1,2,\ldots,s\}$;
\item[(iii)]  $\rk \sf{L_1\mid L_2}{\psi}$ is maximal;
\item[(iv)] $\brk \sf{L_1\mid L_2}{\psi}$ is maximal.
\end{itemize}
\end{lem}

\begin{rem}
Once a basis for every $\CC^{\m(l_i)_t}$ is chosen, we can identify \[\bigoplus_{t=1}^s
\Hom_{\CC}(\CC^{\m(l_1)_t},\CC^{\m(l_2)_t})\] with the space of block-diagonal
matrices $M_{\m(l_1),\m(l_2)}$.
The notation introduced is coherent with Notation \ref{ranknot}. %
\end{rem}

\begin{lem}\label{prod*}
Let $\varphi_1\in (\otimes_{L_1\cup C}W)^G$
and $\varphi_2\in (\otimes_{L_2\cup C}W)^G$ be two tensors. Then,
\begin{itemize}
 \item[(a)] $flat_{L_1\mid L_2}(\varphi_1 \stackrel[C]{}{\ast}  \varphi_2)= flat_{L_1\mid
C}(\varphi_1) flat_{C\mid L_2}(\varphi_2)$;
\item[(b)] $\sf{L_1\mid L_2}{\varphi_1 \stackrel[C]{}{\ast}  \varphi_2}= \sf{L_1\mid C}{\varphi_1}
\sf{C\mid L_2}{\varphi_2}$.
\end{itemize}
\end{lem}
%

%

\begin{notation}
Given an edge $e$, we denote $1_e=\sum_{b\in B} b\otimes b\in \left (W_{e_0}\otimes W_{e_1} \right )^G$. Given a phylogenetic tree $T$, we write $\id_T=(1_e)_{e\in E(T)}$ and call it the \emph{no-mutation presentation of $T$}.
\end{notation}

\subsection{Degenerated trees and trees with observed interior vertices}
For technical reasons, we admit degenerated trees reduced to just
one vertex, which is considered as a leaf.
If $T=\bullet_q$ is such a tree, we associate the $\CC$-vector
space $W$ to $q$, making of $T$ a phylogenetic tree on $(G,W)$. Moreover, we
take $\rep_G(T)$ to be composed of the no-mutation
presentation, that is,
 \begin{eqnarray*}
\rep_G(T)=\{1_q\} \quad \mbox{ where }1_q=\sum_{u\in
B}u\otimes u.
 \end{eqnarray*}
and define  $\Psi_T:\rep_G(T)\ra \LL(T)^G=W$  by mapping $1_q$
to $\sum_{u\in B}u$.
The reader can think of such a tree as a two-leaf tree where we
only accept the no-mutation presentation between its two leaves.
All the above definitions are coherent with this interpretation.

\begin{defn}\label{obser}
Let $T$ be a phylogenetic tree on $(G,W)$ and let $q\in \Int(T)$. Then, we can
write $T=*_{i=1}^m T_i$, where $T_i$ are subtrees of $T$ sharing
the vertex $q$ as a common leaf. Write $T_0$ for the degenerated
tree reduced to $q$. The tree $T$ \emph{with observed} $q$ is
defined by
\[T^q=*_{i=0}^m T_i.\]
Notice that, by definition, the leaves of $T^q$ are the leaves of
$T$,  $L(T^q)=L(T)$, while $\LL(T^q)=\LL(T)\otimes W$. Define a
map \[\Psi_{T}^q:\rep_G(T)\ra \LL(T^q)\] by taking
$\Psi_{T}^q(A)=\Psi_{T_0}(1_q)*\ \Psi_{T_1}(A_1)* \ldots *
\Psi_{T_m}(A_m) $, where $A_i$ is the restriction of the
$G$-evolutionary presentation $A$ to the subtree $T_i$.
\end{defn}
%

\section{The ideal of an equivariant model}
%
In this section, we essentially prove that edge invariants are
indeed phylogenetic invariants (see Introduction). The proof of
this result is quite technical as it is valid for any equivariant
model.

Given a phylogenetic tree $T$ on $(G,W)$ on $W$ and a bipartition $\beta=L_1\mid L_2$
of the leaves of $T$, define $T_1$ (resp. $T_2$) as the minimal
subtree of $T$ that contains the leaves in $L_1$ (resp. $L_2$).
Clearly, we have  $V(T)=V(T_1)\cup V(T_2)$.
Given two vertices $p,q\in V(T)$, the \emph{chain} $ch(q,p)$ is
the linear subtree composed of the edges and vertices between $q$
and $p$.
Define a binary relationship $\sim_{L_1}$ among the leaves of
$L_1$ as follows:
\begin{eqnarray*}
x\sim_{L_1} y \quad \mbox{ if  }\quad ch(x,y)\cap T_2=\emptyset.
\end{eqnarray*}
Analogously, a binary relationship $\sim_{L_2}$ can be defined. It is
easy to check that both $\sim_{L_1}$ and $\sim_{L_2}$ are
equivalence relationships.
Write $n_1$ and $n_2$ for the cardinals of the equivalence classes
of $\sim_{L_1}$ and $\sim_{L_2}$, respectively.
For $i=1,2$, write $\{L_{i,j}\}_{j=1,\ldots,n_i}$ for the
resulting equivalence classes in $L_i$, so that
$L_i=\bigcup_{j=1}^{n_i}L_{i,j}$.
%
Notice that if $l_1=\sharp L_1$ and $l_2=\sharp L_2$, then
$n_1\leq l_1$ and $n_2\leq l_2$.
From now on, we will denote
\[\m_{\beta,T}=\m(\min\{n_1,n_2\}).\]
The main goal of this section is to prove the following
Proposition, which is a generalization of
\cite[Theorem 19.5]{Eriksson05} to equivariant models. Its interest lies in the fact
that it translates the topology of a tree into rank conditions of
suitable matrices.

\begin{prop}\label{rank}
Let $T$ be a trivalent phylogenetic tree $T$ on $(G,W)$ and let $\beta=L_1\mid L_2$ be a bipartition of $L(T)$ as above. Then, we have
\begin{eqnarray*}
\brk \sf{\beta}{\psi}\leq \m_{\beta,T} \qquad \forall \psi\in V(T),
\end{eqnarray*}
and there exists a non-empty Zariski open set $U_{\beta}\subset V(T)$ such that the equality holds for every
$\psi \in U_{\beta}$.
Moreover,
\begin{itemize}
\item[(i)] $\beta$ is an edge  split in $T$ if and only if $\m_{\beta,T}=\m$.
\item[(ii)] If $\beta$ is not an edge  split in $T$, then $\m_{\beta,T} \geq \m(2)$.
\end{itemize}
\end{prop}

The existence of the Zariski open subset above where the
flattening attains the expected rank cannot be proven by a simple
dimension counting as the following example shows.

\begin{example}\rm
Consider $G=\{\rm id\} \subset \mathfrak{S}_4$ and the quartet
tree $T$ having an inner edge $e$. Then $\sf{e}{\psi}$ can be seen
as a $16 \times 16$ matrix $M$ and its expected rank is 4
according to Proposition \ref{rank}(i). The variety $V_G(T)$ has
dimension 60 and is contained in the determinantal variety defined
by the $5\times 5$ minors of $M$, which has dimension
$256-(16-5+1)(16-5+1)=112.$ A priori $V_G(T)$ could also be
included in the variety of $4\times 4$ minors of $M$ which has
dimension $256-(16-4+1)(16-4+1)=87$, so that a general element of
$V_G(T)$ would not have the expected rank 4.
\end{example}

Before proving Proposition \ref{rank}, we need to state a couple of lemmas.

\vspace{2mm}

\begin{lem}\label{trunk}
Let $T$ be a phylogenetic tree on $(G,W)$ and let $\beta=L_1\mid L_2$ be a bipartition of $L(T)$ such that every cherry of $T$ is composed of one leaf in $L_1$ and one leaf in $L_2$. For a generic $G$-evolutionary presentation $A$ of $T$, it holds that $\brk \sf{\beta}{\Psi_T(A)}=\m_{\beta,T}$.
\end{lem}

\begin{proof}
Write $L_1=\{u_1,u_2,\ldots, u_{l_1}\}$ and $L_2=\{v_1,v_2,\ldots, v_{l_2}\}$, and write $n=l_1+l_2$ for the number of leaves of $T$.
Assume that $1\leq l_1\leq l_2$.
Notice that with our assumption, we have  $n_1=l_1$, $n_2=l_2$ and so, $\m_{\beta,T}=\m(l_1)$.
To reach the claim, we first show that the above condition for the rank $\brk$ defined an open set in $\rep_G(T)$. Then, we will prove recursively that this open set is not empty.

Let $\varphi\in \LL(T)^G$ and write
$\sf{\beta}{\varphi}=(\varphi_1,\varphi_2,\ldots,\varphi_s)$.
Then $\sf{\beta}{\varphi}$ has maximal rank if and only
if
\begin{eqnarray}\label{rkB}
\rk \varphi_t=\min \{\m(n_1)_t,\m(n_2)_t\}\quad \mbox{ for  every
}t=1,\ldots, s.
\end{eqnarray}
Each rank condition $\rk \varphi_t<\min \{\m(n_1)_t,\m(n_2)_t\},
t=1,\ldots, s$ defines a closed proper subset $Z_i$ of
$\Hom_G(\otimes_{L_1}W,\otimes_{L_2}W)\cong \LL(T)^G$.
Thus,  \[\widetilde{V}=\LL(T)^G\setminus \cup_{t=1}^s
Z_t\] is a dense open subset in  $\LL(T)^G$, and for every
$\varphi\in \widetilde{V}$,
\[\brk \sf{\beta}{\varphi}= \m_{\beta,T}.\]
Moreover, $\Psi_{T}:\rep_G(T)\ra \LL(T)^G$ is a continuous
map, so $V=\Psi_{T_R}^{-1}(\widetilde{V})$ is an open set in $\rep_G(T)$. To prove that $V$ is non-empty, we will recursively construct a $G$-evolutionary presentation $A\in \rep_G(T)$ with $A_e=1_e$ for any terminal edge $e$, and such that
\begin{eqnarray}\label{equal}
\rk flat_{L_1\mid L_2}\Psi_T(A)=k^{l_1}.
\end{eqnarray}
This implies that the rank of $flat_{L_1\mid L_2}\Psi_T(A)$ is maximal and, by applying (\ref{eq}), we derive that so is $\brk \sf{\beta}{\Psi_T(A)}$. From this, we derive that $\brk \sf{\beta}{\Psi_T(A)}=\m(l_1)=\m_{\beta,T}$.
%

For $n=2$, it is enough to take $A$ equal to the no-mutation presentation: $\id_T=\sum_b b\otimes b$.
For general $n$, take a cherry of $T$.  By reordering the leaves of $L_1$ and $L_2$ if necessary, we can assume that the leaves in this cherry are $u_{l_1},v_{l_2}$. Let $e$ be the edge of $T$ adjacent to it and insert two vertices $q_1$ and $q_2$ in the edge $e$. We obtain a decomposition of $T$ as follows: $T=T^1\ast T^e\ast T^2$, where $T^e$ is a 2-leaf tree with leaves $q_1$ and $q_2$, and $T^1$ and $T^2$ are the subtrees of $T$ obtained when removing $T^e$ from $T$ as shown in figure \ref{dec2}.
Then, we have
\begin{eqnarray*}
L(T^1) & = & \{u_1,u_2,\ldots, u_{l_1-1},q_1,v_1,\ldots,v_{l_2-1}\},\\ L(T^e)& = &\{q_1,q_2\}, \\  L(T^2)& = &\{u_{l_1},v_{l_2},q\}.
\end{eqnarray*}
\begin{figure}
\begin{center}
\psfrag{T}{$T$} \psfrag{T1}{$T^1$} \psfrag{T2}{$T^2$}
\psfrag{Te}{$T^e$}
\psfrag{u1}{\tiny $u_1$}
\psfrag{u2}{\tiny$u_2$}
\psfrag{u3}{\tiny$u_3$}
\psfrag{v1}{\tiny$v_1$}
\psfrag{v2}{\tiny$v_2$}
\psfrag{v3}{\tiny$v_3$}
\psfrag{v4}{\tiny$v_4$}
\psfrag{v5}{\tiny$v_5$}
\psfrag{q1}{\tiny$q_1$} \psfrag{q2}{\tiny$q_2$}
\includegraphics[scale=0.8]{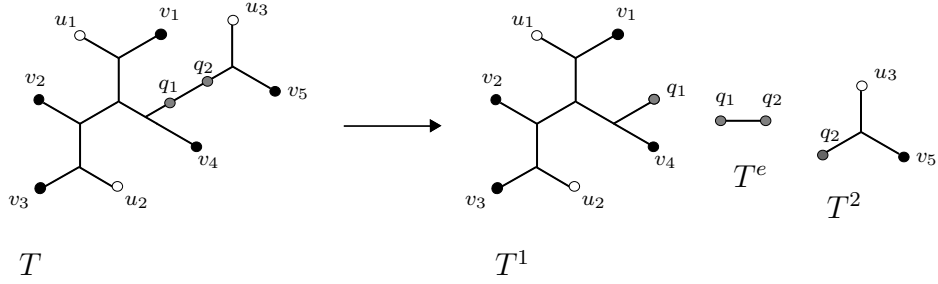}
\end{center}
\caption{\label{dec2} Decomposition of $T$ in the proof of Lemma \ref{trunk}.}
\end{figure}
Write $L^{(1)}_1=\{u_1,\ldots,u_{l_1-1},q_1\}$ and $L^{(1)}_2=\{v_1,\ldots,v_{l_2-1}\}$.
Since $T^1$ has $n-1$ leaves, our assumption says that there is some $G$-evolutionary presentation, say $A^1\in \rep_G(T^1)$, such that
\[\rk flat_{L_1^{(1)} \mid L_2^{(1)}} \Psi_{T^1}(A)=k^{\min  \{ l_1, l_2-1 \}}.\]
Define $A=A^1\ast A_e \ast \id_{T_2}\in \rep_G(T)$, where $A_e\in \Hom_G(W_{q_1},W_{q_2})$ is generic and we will show that the equality (\ref{equal}) holds.
To this aim, we claim that
\begin{eqnarray}\label{claim}
 \Psi_{T}(A) = (\Psi_{T^1}(A^1)\otimes \id)  \stackrel[L_1\cup \{q\}]{}{\ast} \phi(A_e).
\end{eqnarray}
\emph{Proof.}
First of all, the decomposition $A=A^1\ast A_e \ast \id_{T_2}$ induces a decomposition of
$\Psi_{T}(A)$ as
\begin{eqnarray*}\label{tut}
\nonumber \Psi_{T}(A) & = & \Psi_{T^1}(A^1)\ast \Psi_{T^e}(A_e) \ast \Psi_{T^2}(\id_{T_2})=\\
\nonumber & = & \sum_{z_1,z_2\in B} (A_e\mid z_1\otimes z_2) \langle\Psi_{T^1}(A^1)\mid z_1\rangle \otimes \langle \Psi_{T^2}(\id_{T_2})\mid z_2\rangle.
\end{eqnarray*}
and notice that
\begin{eqnarray*}
\lefteqn{ \langle \Psi_{T^1}(A^1) \mid z_1 \rangle \otimes \langle \Psi_{T^2}(\id_{T_2}) \mid z_2 \rangle = \langle \Psi_{T^1}(A^1) \mid z_1 \rangle \otimes  (z_2 \otimes z_{2} ) =} \\ & &
\sum_{\substack{b_i\in B_{u_i}\\i=1,\ldots,l_1-1} }\langle \Psi_{T^1}(A^1)\mid b_1\otimes \ldots \otimes b_{l_1-1}\otimes z_1\rangle \otimes ( z_2\otimes z_2) \otimes b_1\otimes \ldots \otimes b_{l_1-1}.
\end{eqnarray*}
\normalsize
Thus, we obtain that $\Psi_T(A)$ is equal to
\small
\begin{eqnarray*}
 \sum_{z_1,z_2,b_1,\ldots,b_{l_1-1}} (A_e\mid z_1\otimes z_2)\langle \Psi_{T^1}(A^1)\mid b_1\otimes \ldots \otimes b_{l_1-1}\otimes z_1\rangle \otimes z_2 \otimes z_2 \otimes b_1\otimes \ldots \otimes b_{l_1-1} .
\end{eqnarray*}
\normalsize
On the other hand, consider the $G$-equivariant map  
\begin{eqnarray*}
\varphi(A_{e}): \otimes_{L_1\cup \{q_1\}}W & \ra &  \otimes_{L_1} W \\
b_1\otimes \ldots \otimes b_{l_1-1}\otimes b_{l_1} \otimes b_{q_1} & \mapsto &
( A_e\mid b_{l_1} \otimes  b_{q_1}) b_1\otimes \ldots \otimes b_{l_1-1}\otimes  b_{q_1}.
\end{eqnarray*}
and define the tensor $\phi(A_e)$ as the image of $\varphi(A_e)$ via the isomorphism \[\Hom_{G}(\otimes_{L_1\cup \{q_1\}}W,\otimes_{L_1} W)\cong \left ( ( \otimes_{L_1\cup \{q_1\}}W ) \otimes ( \otimes_{L_1} W ) \right )^G.\]
Thus, if $b_i\in B_{u_i}, i=1,\ldots, l_1-1$,  $z_1\in B_{q_1}$ and $z_2\in B_{l_1}$, then
\begin{eqnarray*}
 \langle \phi(A_e)\mid b_1\otimes \ldots \otimes b_{l_1-1}\otimes z_1\otimes z_2\rangle =( A_e\mid z_1 \otimes z_2 ) z_2 \otimes b_1\otimes \ldots \otimes b_{l_1-1}.
\end{eqnarray*}
If $\id=\sum_b b\otimes b\in W_{u_{l_1}}\otimes W_{u_{l_2}}$, we have $\langle \id \mid z_2 \rangle = z_2$ and so 
\begin{eqnarray*}
 \langle \Psi_{T^1}(A^1)\otimes \id \mid b_1\otimes \ldots \otimes b_{l_1-1}\otimes z_1\otimes z_2\rangle
= \langle \Psi_{T^1}(A^1)\mid b_1\otimes \ldots \otimes b_{l_1-1} \otimes z_1 \rangle \otimes z_2.
\end{eqnarray*}
Putting all together, we obtain that
$(\Psi_{T^1}(A^1)\otimes \id)  \stackrel[L_1\cup \{q\}]{}{\ast}
\phi(A_e)$ is equal to \small
\begin{eqnarray*}
\sum_{b_1,\ldots,b_{l_1-1},z_1,z_2} (A_e\mid z_1\otimes z_2 )\langle \Psi_{T^1}(A^1) \mid b_1\otimes \ldots \otimes b_{l_1-1} \otimes z_1 \rangle \otimes z_2 \otimes \left ( z_2  \otimes  b_1\otimes \ldots \otimes b_{l_1-1} \right ).
\end{eqnarray*}
\normalsize
This proves the claim. 

Now, apply Lemma \ref{prod*} to (\ref{claim}) to get
\[flat_{L_1\mid L_2}\Psi_{T}(A) = flat_{L_1\mid L_1\cup \{q\}}(\Psi_{T^1}(A^1)\otimes \id) flat_{L_1\cup \{q\}\mid L_2} \phi(A_e).\]
%
%
It is straightforward to check that
\[flat_{L_1\cup \{q\} \mid L_2}\left (\Psi_{T^1}(A^1)\otimes \id \right )=\left (flat_{L^{(1)}_1 \mid L^{(1)}_2}\Psi_{T^1}(A^1) \right ) \otimes \left (flat_{\{u_{l_1}\}\mid \{v_{l_2}\}}\id \right),\]so the rank of $F:=flat_{L_1\cup \{q\} \mid L_2}\left (\Psi_{T^1}(A^1)\otimes \id \right )$ is equal to the product of ranks:
$k \times k^{\min\{l_1-1,l_2\}}=k^{\min\{l_1,l_2+1\}}
$.
On the other hand, write $G(A_e)$ for the matrix of $flat_{L_1\mid L_1\cup \{q\}}\phi(A_e)$ in the basis $B(\otimes_{L_1}W)$ and $B(\otimes_{L_1\setminus \{q\}}W)$.
It is a block diagonal matrix, each block being a convenient column of the matrix $A_e$.
Then the  rank of $flat_{L_1\mid L_2}\Psi_T(A)$ is $k^{l_1}$ if and only if
$\ker(F)\cap \im(G(A_e)) =\{0\}$.
Since this holds for a generic matrix $A_e$, the claim follows.
\end{proof}

\begin{lem}\label{observed}
Let $T$ be a phylogenetic tree on $(G,W)$ and let $q\in \Int(T)$. Assume that
$q$ has degree two while the remaining interior vertices have
degree three and write $T^q$ for the tree with observed $q$. Then, for a generic $G$-evolutionary presentation $A\in \rep_G(T)$, we have
that
\begin{eqnarray*}
\brk \sf{\{q\}\mid L(T)}{\Psi_{T^q}(A)}=\m.
\end{eqnarray*}
\end{lem}

\begin{proof}
First of all, notice that $\brk \sf{\{q\}\mid L(T)}{\psi}\leq \m$
for all $\psi \in \LL(T^q)^G$ and that there is a non-empty open set
$U\subset \LL(T^q)^G$ where the above
inequality holds.
Indeed, if $\sf{\{q\}\mid L(T)}{\psi}=(\psi_1,\ldots,\psi_s)$ it is
enough to take $U=\cap_{i=1}^s U_i$, where each $U_i$ is defined
by asking that $\rk \psi_i=m_i$. Every $U_i$ is a dense open
set in $\LL(T)^G$, and so is $U$.

To reach the claim we only have to prove that $\Psi_{T^q}^{-1}(U)$
is not empty. To this aim, it will be enough to consider the
no-mutation presentation $\id=\{1_e\}_{e\in E(T^q)}$.
The linear map $flat_{\{q\}\mid L(T)}\Psi_{T^q}(\id):W\ra
\otimes_{L(T)} W$ defined by $b\mapsto  b\otimes
\ldots  \otimes b$ has rank equal to $\dim(W)$. By
virtue of (\ref{eq}), we infer that $\id\in \Psi_{T^q}^{-1}(U)$
and we are done.
%
%
\end{proof}

Now, we come back to the general case. So, let $T$ be a
phylogenetic tree on $(G,W)$ and let $\beta=L_1\mid L_2$ be a
bipartition of $L(T)$. We introduce some terminology and notation
that will be helpful. For the seek of clarity, this notation will
not reflect its dependence on $\beta$, but confusion should not
arise since the bipartition is fixed throughout this section. Keep
the notation introduced at the beginning of this section:

\begin{figure}
\begin{center}
\psfrag{T}{$T$} \psfrag{u11}{\tiny $u_1^1$}
\psfrag{u12}{\tiny$u_1^2$} \psfrag{u21}{\tiny$u_2^1$}
\psfrag{u22}{\tiny$u_2^2$} \psfrag{u31}{\tiny$u_3^1$}
\psfrag{A}{\tiny$T_R$} \psfrag{T11}{\tiny$T_1^1$}
\psfrag{T12}{\tiny$T_1^2$} \psfrag{T21}{\tiny$T_2^1$}
\psfrag{T31}{\tiny$T_3^1$} \psfrag{T22}{\tiny$T_2^2$}
\includegraphics[scale=0.8]{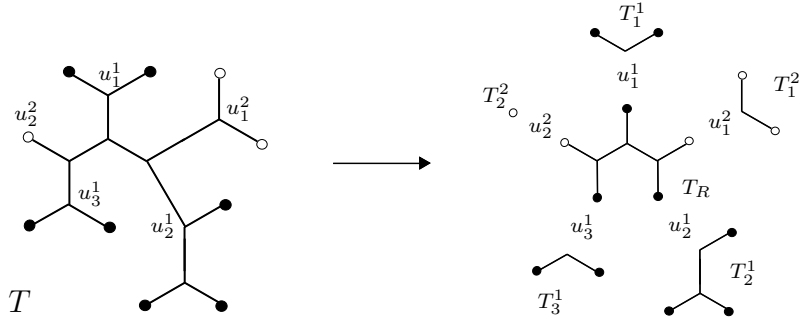}
\end{center}
\caption{\label{trees} Construction of the trunk and the boughs
of a given tree $T$. Black and white dots represent the leaves in
$L_1$ and $L_2$, respectively. Notice that some $T_{i,j}$ may be
reduced to the vertex $u_{i,j}$. Indeed, this happens if and only
if $u_{i,j}\in L_i$.}
\end{figure}

\begin{itemize}
 \item[1.]For $i=1,2$ and $j=1,\ldots, n_i$, denote by $T_{i,j}$ the minimal phylogenetic subtree of $T$
containing $L_{i,j}$. These subtrees are
called the \emph{boughs of $T$ relative to $\beta$}.
Every $T_{i,j}$ has a distinguished vertex, denoted by $u_{i,j}$,
which has degree two. All the remaining interior vertices have
degree 3.
For $i=1,2$, write $L_i^R=\{u_{i,j}\}_{j=1,\ldots,n_i}$.

\item[2.] The \emph{trunk of $T$ relative to $\beta$}, denoted by $T_R$, is the phylogenetic tree on $(G,W)$
obtained when removing all the boughs from $T$. Equivalently,
$T_R$ is the minimal subtree of $T$  containing all the $u_{i,j}$,
so that $L(T_R)=L_1^R\cup L_2^R$.

 \end{itemize}
See figure \ref{trees} for an example of this construction.

Notice that $T_R$ is reduced to a 2-leaf tree if and only if
$\beta$ is an edge split of $T$.
Notice also that every cherry in $T_R$ is composed of one leaf in
$L_1^R$ and one leaf in $L_2^R$.

\vspace{2mm}
\begin{proofof}{Proposition \ref{rank}}
Every $G$-evolutionary presentation ${A\in
\rep_G(T)}$ induces by restriction a $G$-evolutionary presentation
$A^R$ in the trunk $T_R$ and $G$-evolutionary presentations $A_{i,j}$
in the boughs $T_{i,j}$. Actually, the mappings $A\mapsto A^R$
and $A\mapsto A_{i,j}$ define continuous maps
\begin{eqnarray*}\label{AR}
\pi_R:\rep_G(T) \rightarrow \rep_G(T_R) \qquad \pi_{i,j}:\rep_G(T)
\rightarrow \rep_G(T_{i,j}).
\end{eqnarray*}

\vspace{3mm} Given $A\in \rep_G(T)$, we proceed to decompose
$\Psi_T(A)$ in terms of tensors associated to the trunk and the
boughs of $T$.
To this aim and for every $i,j$, consider the tree
$T_{i,j}^{u_{i,j}}$ with the vertex $u_{i,j}$ observed in it (see
section 2.2). It is straightforward to check that $T$ is recovered by joining these boughs to the trunk:
\[T=T_{1,1}^{u_{1,1}} \ast ( \ldots \ast ( T_{1,n_1}^{u_{1,n_1}} \ast (T_{2,1}^{u_{2,1}}\ast  (\ldots \ast (T_{2,n_2}^{u_{2,n_2}}\ast T_R )\ldots )) \ldots ).\]
Regarding $A_{i,j}$ as a $G$-evolutionary presentation of
$T_{i,j}^{u_{i,j}}$, write $\varphi_{i,j}(A)\in \LL(T_{i,j}^{u_{i,j}})$
for the image of $A_{i,j}$ by the map $\Psi_{T_{i,j}}^{u_{i,j}}$
defined in Definition \ref{obser}. Then write
\begin{eqnarray*}
\varphi_1^A& = &\otimes_{j=1}^{n_1} \varphi_{1,j}(A)\in   (\otimes_{L^R_1\cup L_1} W)^G \\
\varphi_2^A& = &\otimes_{j=1}^{n_2} \varphi_{2,j}(A)\in  (\otimes_{L^R_2\cup L_2} W)^G\\
\varphi_R^A & = & \Psi_{T_R}(A^R)\in \LL(T_R)^G.
\end{eqnarray*}
From the construction of these three tensors, it is clear that
(see Notation \ref{new}):
\[\Psi_T(A)=\varphi_1^A*\varphi_{R}^A*\varphi_2^A.\]
Write $\beta_1=L_1\mid L_1^R$, $\beta_2=L_2^R \mid L_2$ and $\beta_R=L_1^R\mid L_2^R$.
By applying (\ref{prod*}) we infer the following equality of
maps
\begin{eqnarray*}
\label{equal_matrices} \sf{\beta}{\Psi_T(A)}=
\sf{\beta_1}{\varphi_{1}^A}   \sf{\beta_R}{\varphi_R^A}  \,  \sf{\beta_2}{\varphi_2^A}
\end{eqnarray*}
and from it, the inequality $\brk \sf{\beta}{\Psi_T(A)}\leq
\m_{\beta,T}$.
Next, we prove that the equality  actually holds for a generic
$G$-evolutionary presentation.
To this aim, we will show that the three maps above have
maximal rank for a generic
$G$-evolutionary presentation $A$ of $T$. It
is straightforward to check that this will imply that $\brk
\sf{\beta}{\Psi_T(A)}=\m_{\beta,T}$  (use for instance
the Frobenius inequality, see Section 2.9.6 of \cite{EMT}).
%


First of all, we infer from Lemma \ref{trunk} that there is a dense open set $V_R\subset \rep_G(T_R)$ such that for every $B\in V_R$, $\brk \sf{\beta_R}{\Psi_{T_R}(B)}=\m_{\beta_R,T_R}$.
Since the map $\pi_R:\rep_G(T)\ra \rep_G(T_R)$ is surjective, the
set \[\mathbf{U}_R:=\pi_R^{-1}(V)\]is a \emph{non-empty} Zariski
open set in $\rep_G(T)$.

Similarly, for every $i=1,2$ and $j=1,2,\ldots, n_i$, let $V_{i,j}\subset
\rep_G(T_{i,j})$ be the dense open set defined by Lemma \ref{observed}
applied to $T_{i,j}^{u_{i,j}}$ and $U_{i,j}=\pi_{i,j}^{-1}(V_{i,j})$. Since $\rep_G(T)$ is irreducible, it is clear that
\[\mathbf{U}_i=\bigcap_{j=1}^{n_i}U_{i,j},\quad i=1,2\]is non-empty and open.
 On the other hand, for any $A\in
\rep_G(T)$, we have
\[\sf{\beta_i}{\varphi_{i,1}(A)\otimes \ldots \otimes  \varphi_{i, n_i}(A)}= \otimes_{j=1}^{n_i}  \sf{L_{i,j}\mid \{u_{i,j}\}} {\varphi_{i,j}(A)} \] and therefore (see for instance \cite{EMT})
\begin{eqnarray*}
\rk \,\sf{\beta_i}{\varphi_i^A}=\prod_{j=1}^{n_i} \rk \,
\sf{L_{i,j}\mid \{u_{i,j}\}}{\varphi_{i,j}(A)}=\m(n_i).
\end{eqnarray*}
Thus, if $A\in \mathbf{U}_i$,  the rank of
$\sf{L_1,L^R_1}{\Psi_1^A}$ is maximal.
%

To finish the proof it is enough to take 
\[U_{\beta}=\Psi_{T}(\mathbf{U}_1 \cap \mathbf{U}_R \cap \mathbf{U}_2)\subset V_T.\]
Finally, if $\beta=L_1\mid L_2$ is an edge split, then $n_1=n_2=1$ and $T_R$ is 2-leaf
tree. It follows that $\m_{\beta,T}=\m$. This proves (i). If
$\beta$ is not an edge split, it is clear that  $n_1,n_2\geq
2$ and the claim of (ii) follows by Lemma \ref{max_m}.
\end{proofof}

\begin{rem}\label{stochastic}
The preceding proof actually shows that the dense open set
$U_{L_1\mid L_2}\subset \rep_G(T)$ cuts the set of stochastic
parameters, i.e
\[U_{L_1\mid L_2} \cap \prod_{e\in E(T)} \Delta^G\neq \emptyset,\]where $\Delta^{G}$ is the set of Markov matrices, that is, matrices whose entries are all non-negative and whose columns sum to 1.
Indeed, as suggested by Lemma \ref{trunk} and Lemma
\ref{observed}, it is enough to take $A\in \rep_G(T)$ with
$A_e=\sum_{b\in B} b\otimes b$ whenever $e\in E(T_R)$ is terminal
or $e\in E(T_{i,j})$ for some $i,j$, and $A_e$ a generic Markov matrix,
otherwise.
\end{rem}

Proposition \ref{rank} suggests the following definitions.

\begin{defn}
If $L_1\mid L_2$ is a bipartition of $L(T)$, the \emph{ideal of}
$L_1\mid L_2$, denoted by $I_{L_1 \mid L_2}$, is the ideal in the
coordinate ring of $\LL(T)^G$ defined by the conditions
\[\brk \sf{L_1\mid L_2}\psi \leq \m\]being $\psi \in \LL(T)^G$  a tensor of indeterminates. Equivalently, $I_{L_1 \mid L_2}$ is generated  by the
$(m_t+1)$-minors of the $t$-th box of $\sf{L_1 \mid L_2}\psi \in
M_{\m(l_1), \m(l_2)}$, for  $t=1, \dots,s$ (see
Notation \ref{basis/flat}).
\end{defn}

\begin{notation}
Let $T$ be a phylogenetic tree on $(G,W)$ and let $e$ be an edge of $T$
that splits the leaves into two sets $L_1$ and $L_2$ of
cardinality $l_1$ and $l_2$, respectively. The ideal  $I_{L_1\mid
L_2}$ will be also denoted as $I_e$. Due to Proposition \ref{rank}
we have that if $e$ belongs to $E(T)$, then $I_e\subseteq \I(T)$.
\end{notation}

\begin{defn}
The \emph{edge invariants} of $T$  are the elements of the ideal
$\sum_{e\in E(T)}I_e$.
\end{defn}

Proposition \ref{rank} proves that edge invariants are
\textit{phylogenetic invariants}, that is, elements in
$\mathcal{I}(T)$  that do not vanish on all points of $\cup_{T}
V(T)$ where the union runs over all trivalent tree topologies.
Indeed, given a G-spaced tree $T_0$ and an edge $e \in E(T_0)$,
there exist trivalent trees that do not have $e$ as an edge split
and so $I_e$ is not contained in $\mathcal{I}(\cup_{T} V(T))$.

Is is worth highlighting that using Proposition \ref{rank} we also obtain the \textit{generic identifiability} of the tree topology for equivariant models. The tree topology of a model of sequence mutation is said to be generically identifiable if for generic choices of stochastic parameters $A \in \prod_{e\in E(T)} \Delta^G, A' \in \prod_{e\in E(T')} \Delta^G$ (see Remark \ref{stochastic}), $\Psi_T(A)=\Psi_{T'}(A')$ implies $T=T'$ (see for instance \cite{Allman2006}).  In order to prove this kind of results, one only has to show the corresponding irreducible varieties $V(T)$ and $V(T')$ are not contained one into the other.%
We obtain the following result that was already known for the
general Markov model (see \cite{Steel94}) and for group-based
models \cite{Steel1992}.

\begin{cor}\label{identifiability}
The tree topology is generically identifiable in all equivariant
evolutionary models.
\end{cor}

\begin{proof}
Let $T,T'$ be two different trivalent phylogenetic tree on
$(G,W)$. Then there is an edge split $e$ in $T$ that is not an
edge split in $T'$. By Proposition \ref{rank}, there exists an
element $f$ in $I_e$ (and therefore in $\mathcal{I}(T)$) that does
not belong to $\mathcal{I}(T')$.  In terms of varieties this
proves that $V(T') \subsetneq V(T)$, and that $V(T) \subsetneq
V(T')$ is proven similarly. As $V(T)$ and $V(T')$ are irreducible
varieties, this shows that they meet properly.
\end{proof}

\section{Phylogenetic Invariants}

The purpose of this section is to prove that, for phylogenetic
reconstruction, the only relevant invariants are the edge
invariants introduced in the previous section. This is a natural
result if one takes into account the Splits Equivalence Theorem in
combinatorics (see Theorem \ref{splits} below). Let $\mathcal{T}$
be the set of trivalent tree topologies with leaf set $\{v_1,v_2,
\dots, v_n\}$. Two bipartitions $L_1|L_2$, $M_1|M_2$ of a  set $L$
are said to be \textit{compatible} if at least one of the four
intersections $L_1 \cap M_1$, $L_1\cap M_2$, $L_2\cap M_1$,
$L_2\cap M_2$ is empty. For example, if $L_1|L_2$, $M_1|M_2$ are
two edge splits of the same tree $T$, then they are compatible. We
recall that any trivalent tree on $n$ leaves has $2n-3$ interior
edges.

\begin{thm}[\cite{Buneman}, {\cite[Theorem 2.35]{ASCB2005}}]
\label{splits}
A collection $\mathcal{B}$ of $2n-3$ bipartitions is pairwise compatible if and
only if there exists a tree $T \in \mathcal{T}$ such that $\mathcal{B}$ is the set of edge
splits of $T$. Moreover, if such a tree $T$ exists then it is
unique.
\end{thm}

In order to make our result concerning phylogenetic invariants
more precise we need to introduce some notation.

We fix $G \subset \mathfrak{S} _k$ and $W$ as in section 2 and each topology $T \in \mathcal{T}$ will be considered as a phylogenetic tree on $(G,W)$. Then all trees $T$ in $\mathcal{T}$ have the same space of $G$-tensors 
which will be denoted by $\mathcal{L}=(\bigotimes_{i=1}^{n} W)^G$.

\begin{defn}
Let  $\o$ be an $s$-tuple and let $\beta=L_1\mid L_2$ be a
bipartition of $\{v_1,v_2, \dots, v_n\}$. Then we let $D_{\leq
{\o}}^{\beta}$ be the subvariety of $\LL$ defined as
$$D_{\leq \mathbf{o}}^{\beta}=\{\psi \in \LL \mid \brk \sf{\beta}\psi \leq \mathbf{o} \}\,$$
and, if the thin flattening of $\psi \in \mathcal{L}$ is
$\sf{\beta} \psi=(\psi_1,\psi_2,\ldots,\psi_s)$, we define
$D_{< {\o}}^{\beta}$ to be the set $$D_{< {\o}}^{\beta}= \{\psi
\in \mathcal{L} \mid \rk \psi_j < o_j \textrm{ for some j } \}.
$$  For example, $D_{\leq \mathbf{m}}^{\beta}$ coincides with the
set of zeroes $Z(I_{L_1,L_2})$. Notice that both $D_{\leq
\mathbf{o}}^{\beta}$ and $D_{< \mathbf{o}}^{\beta}$ are algebraic
sets although the second is not irreducible.
%
\end{defn}

%
%


\begin{notation}
Given a tree $T\in \mathcal{T}$ and using the notation of Proposition \ref{rank},
for each bipartition $\beta=L_1 \mid L_2$ of $\{v_1,v_2, \dots, v_n\}$,
we call $\m_{\beta,T}$ 
the maximum rank that $\sf{\beta}{\psi}$ can have if $\psi$
belongs to $V(T)$.
Then Proposition \ref{rank} shows that  
$$V(T) \subseteq D^{\beta}_{\leq \m_{\beta,T}}$$
and that $V(T) \setminus D^{\beta}_{< \m_{\beta,T}}$ is a dense
open subset of $V(T)$ for any bipartition $\beta=L_1 \mid L_2$. We
call this open subset $U_{T,\beta}$, so that $U_{T,\beta}=V(T)
\setminus D^{\beta}_{< \m_{\beta,T}}$ is the locus of tensors
$\psi \in V(T)$ that satisfy $\brk
\sf{\beta}{\psi}=\m_{\beta,T}$. We define
$U_T=\cap_{\beta}U_{T,\beta}$, where the intersection is taken
among all bipartitions of $\{v_1,v_2, \dots, v_n\}$. As $V(T)$ is
an irreducible variety, $U_T$ is still a dense open subset of
$V(T)$ and it corresponds to the set of points in $V(T)$ whose
flattening $\sf{\beta}{\psi}$ along any partition $\beta$ of the
set of leaves of $T$ has the expected rank $\m_{\beta,T}$.
\end{notation}

With this set up in mind, the main result of this paper is the
following.

\begin{thm}\label{teophyloinvar} For each $T \in
\mathcal{T}$, let $U_T \subset V(T)$ be the dense open set defined
above. Let $p$ be a point in $\bigcup_{T \in \mathcal{T}}U_T
\subseteq \LL$ and let $T_0$ be any tree in $\mathcal{T}$. Then,
$p$ belongs to $V(T_0)$ if and only if $p$ belongs to the set of
zeroes $Z(\sum_{e \in E(T_0)}I_e)$.
\end{thm}

\begin{rem}
As we pointed out in the introduction, this result says that for a
general point on $\bigcup_{T\in \T}V(T)$, it is enough to evaluate
the edge invariants to decide to which variety $V(T)$ the point
actually belongs to.

This result would still hold for non-trivalent trees when imposing
that all trees in the corresponding set $\mathcal{T}$ have the
same collection of degrees at interior vertices.
\end{rem}

After all the technical issues in section 3, the proof of Theorem
\ref{teophyloinvar} is now straightforward.

\begin{proofof}{\ref{teophyloinvar}}
By Proposition \ref{rank} we already know that $\sum_{e \in
E(T_0)}I_e \subseteq \I(T_0)$, therefore if $p \in V(T_0)$, we
immediately have that $p$ belongs to $Z(\sum_{e \in E(T_0)}I_e)$.

Conversely, let $p\in \cup_{T \in \mathcal{T}} U_T.$ Then $p$
belongs to $U_T \subset V(T)$ for a certain $T \in \mathcal{T}$,
so that $\brk \sf{\beta} p = \m_{\beta,T}$ for any bipartition
$\beta$ of $\{v_1,v_2, \dots, v_n\}$. On the other hand, if $p \in
Z(\sum_{e \in E(T_0)}I_e)$, then $p \in Z(I_e)$ for any $e \in
E(T_0)$ and hence, $\brk \sf{e}{p} \leq \m$ for all $e \in E(T_0)$.
This implies that $\m_{e,T} \leq \m$ for all $e \in E(T_0)$, which
can only happen if $e$ is a split of $T$ for all $e \in E(T_0)$
(see Proposition \ref{rank}). But two trivalent trees $T$ and
$T_0$ on $n$ leaves have the same collection of splits if and only
if $T=T_0$ (see Theorem \ref{splits}), so the proof is concluded.
\end{proofof}

\begin{rem} Theorem \ref{teophyloinvar} also says that the intersection $U_T \cap U_{T'}$ is empty for any $T \neq T' \in \mathcal{T}$.
However, there exists points  in $V(T) \cap V(T')$ for any $T \neq T'$. Indeed, it is enough to consider $\psi_T(A)$ where $A$ is the no-mutation presentation; then $\psi_T(A)$ lies in $V(T')$ for all $T'$. This proves that $\bigcap_T V(T)$ is not empty but one can also prove that, if $n \geq 5$, for any two different tree topologies $T_1, T_2$ one has $V(T_1) \cap V(T_2) \neq \bigcap_T V(T).$
 \end{rem}

In the next Corollary we give an open subset $\mathcal{U}$ defined
intrinsically from the ambient space $\mathcal{L}$ such that
$\mathcal{U} \cap \cup_T V(T) = \cup_T U_T$. This is
relevant for biological applications because then we will be able
to check whether the given data point lies (or rather \textit{is
close to}) in $\cup_T U_T$. From now on let $\mathcal{B}$ be the set of all bipartitions of $\{v_1, \dots,v_n\}.$

\begin{cor}\label{corphyloinvar}
Let $\mathcal{U}=\bigcup\limits_{T \in \mathcal{T}} \bigcap\limits_{\beta \in \mathcal{B}} (\LL \setminus
D^{\beta}_{<\m_{\beta,T}}).$
Then
$$\mathcal{U} \cap \bigcup\limits_{T \in \mathcal{T}}V(T) = \bigcup_{T \in \mathcal{T}} U_T$$
and
if  $p$ is a point in $\mathcal{U} \cap \bigcup_{T \in
\mathcal{T}}V(T)$ and $T_0$ is any tree in $\mathcal{T}$, then $p$
belongs to $V(T_0)$ if and only if $p$ belongs to the set of
zeroes $Z(\sum_{e \in E(T_0)}I_e)$.
\end{cor}

\begin{proof}
We just need to prove that $\mathcal{U}\cap(\bigcup_{T \in
\mathcal{T}}V(T)) = \bigcup_{T \in \mathcal{T}}U_T$ because the other
assertion follows from Theorem \ref{teophyloinvar}.

We have $\mathcal{U}\cap (\bigcup_{T \in
\mathcal{T}}V(T))=\bigcup_{T,T'} V(T) \cap (\cap_{\beta} \LL
\setminus D^{\beta}_{< \m_{\beta,T'}}) $. If $T \neq T'$ this
intersection is the empty as we can see taking $\beta$ an edge
split of $T$ but not of $T'$. Hence we obtain $\mathcal{U} \cap
(\cup_{T \in \mathcal{T}}V(T))=\bigcup_{T} V(T) \cap
(\bigcap_{\beta} \LL \setminus D^{\beta}_{< \m_{\beta,T}}) $,
which is precisely $\cup_T U_T$.
\end{proof}

In terms of ideals, Theorem \ref{teophyloinvar} says the
following:
\begin{cor}\label{corideals} Let $R$ be the polynomial ring of $\LL$ and let $f$ be
any element
in \[\left(\sum\limits_{T \in \mathcal{T}} \bigcap\limits_{\beta \in \mathcal{B}}\I(D_{<\m_{\beta,T}}) \right) \setminus \bigcap\limits_T\I(T).\]
Then, the following equality holds in the localized ring $\left(R\diagup \bigcap_T \I(T)\right)_{\overline{f}}$
$$\bigg( \I(T_0)\diagup \bigcap_T\I(T)\bigg)_{\overline{f}} =\left(rad\bigg(\sum_{e \in E(T_0)}I_e\bigg) \diagup \bigcap_T \I(T)\right)_{\overline{f}}.$$
\end{cor}
\begin{proof}
If we are given an $f$ as above, then $U_f:= \LL \setminus
\{f=0\}$ is contained inside the open set $\mathcal{U}$ defined in Corollary
\ref{corphyloinvar}. Indeed, an $f$ as above is contained inside $rad(\sum_{T \in \mathcal{T}} \cap_{\beta}\I(D_{<\m_{\beta,T}}))$ which is equal to $\I(\cap_{T} \cup _{\beta} D^{\beta}_{<\m_{\beta,T}}).$ Therefore $\cap_{T} \cup _{\beta} D^{\beta}_{<\m_{\beta,T}} \subset \{f=0\}$ and $U_f \subset \LL \setminus \cap_{T} \cup _{\beta} D^{\beta}_{<\m_{\beta,T}}=\mathcal{U}$.


In particular, $U_f \cap (\cup_T V(T))$ is contained inside $\cup_T U_T$. Therefore in $U_f$ we still have that the variety $V(T_0)$ is
defined inside $\cup_{T \in \mathcal{T}}V(T)$ by $\sum_{e \in
E(T_0)}I_e$. Hence in terms of ideals in $R_f$ we obtain the
equality above.
\end{proof}

We do not know whether $\sum_{e \in E(T_0)}I_e$ is a radical ideal
so we cannot remove $rad$ from the expression above. We pose the following question:

\vspace{2mm}

\begin{quest}
Given a set $S$ of compatible splits, is $\sum_{\beta\in S}
I_{\beta}$ radical?
\end{quest}

\begin{rem}
In order to check whether Theorem \ref{teophyloinvar} can be
applied to a given data point $p \in \LL$, it is enough to check
that $f(p) \neq 0$ for a generic $f$ in \[\left(\sum\limits_{T \in
\mathcal{T}} \bigcap\limits_{\beta \in
\mathcal{B}}\I(D_{<\m_{\beta,T}}) \right) \setminus
\bigcap\limits_T\I(T).\] Such a polynomial $f$ should be chosen a
priori, so that when dealing with data one does not need to
compute this ideal.
\end{rem}

\begin{rem}
It is interesting to explore whether $U_T$ can be defined by a complete intersection in the sense of \cite{CFS2}. This would reduce the number of generators of $I_e$ to be used in phylogenetic reconstruction. However, this is another issue on which we plan to work in the future.
\end{rem}

Although the degrees of a set of generators of the ideal of a
phylogenetic tree evolving under the general Markov model or under
the strand symmetric model are not known, Theorem
\ref{teophyloinvar} allows us to give the degrees of those
invariants that are relevant in phylogenetic reconstruction. It is
worth highlighting that these degrees do not depend on the number
of leaves but only on the model and can be computed a priori (see
the next sections for the precise examples of evolutionary
models).

\begin{cor}\label{degrees} Let $(G, W)$ be an equivariant evolutionary model and let $\m=(m_1, \dots, m_s)$ be defined as in section 3. Then, for any tree topology on any number of leaves, the polynomials that are relevant for recovering the tree topology in phylogenetics have degrees in  $\{m_1+1, \dots, m_s+1\}$. In particular, the relevant phylogenetic invariants for the following evolutionary models have degrees:
\begin{itemize}
\item 5 for the general Markov model.
\item 3 for the strand symmetric model.
\item 2 for the Kimura 3-parameter model.
\item 1 or 2 for the Kimura 2-parameter model.
\item 1 or 2 for the Jukes-Cantor model.
\end{itemize}
\end{cor}

\section{Examples}

In this section, we study some well-known evolutionary models in
phylogenetics. Let $B=\{\a,\c,\g,\t\}$ be  the set of the four
nucleotides and take $W=\langle \a,\c,\g,\t \rangle_{\CC}\cong
\CC^4$ with the bilinear form $(\cdot \mid \cdot)_W$ that makes
$B$ orthonormal.
We consider the group of permutations of 4 elements,
\[\mathfrak{S}_4=Sym \{B\}.\]It is generated by $g_1=(\textrm{id}),
g_2=(\a\c), g_3=(\a\c\g), g_4=(\a\c\g\t)$ and  $g_5=(\a\c)(\g\t)$,
which correspond to the five conjugacy classes of
$\mathfrak{S}_4$. We work with the natural permutation linear
representation  ${\rho:\mathfrak{S}_4 \rightarrow \rm{GL}(W)}$
given by permuting the coordinates of~$W$:
\begin{eqnarray*}
\tiny g_1 \mapsto \left ( \begin{array}{cccc} 1 & 0 & 0 & 0 \\ 0 &
1  & 0 & 0 \\ 0 & 0 & 1 & 0 \\ 0 & 0 & 0  & 1 \end{array} \right )
& \tiny g_2 \mapsto \left ( \begin{array}{cccc} 0 & 1 & 0 & 0 \\ 1
& 0  & 0 & 0 \\ 0 & 0 & 1 & 0 \\ 0 & 0 & 0  & 1 \end{array} \right
) &
\tiny g_3\mapsto  \left ( \begin{array}{cccc} 0 & 0 & 1 & 0 \\ 1 & 0  & 0 & 0 \\ 0 & 1 & 0 & 0 \\ 0 & 0 & 0  & 1 \end{array} \right ) \\
\tiny g_4 \mapsto \left ( \begin{array}{cccc} 0 & 0 & 0 & 1 \\ 1 &
0  & 0 & 0 \\ 0 & 1 & 0 & 0 \\ 0 & 0 & 1  & 0 \end{array} \right )
& \tiny g_5 \mapsto \left ( \begin{array}{cccc} 0 & 1 & 0 & 0 \\ 1
& 0  & 0 & 0 \\ 0 & 0 & 0 & 1 \\ 0 & 0 & 1  & 0  \end{array}
\right )
\end{eqnarray*}
Write $\chi=\rm{Tr}(\rho(\cdot))$ for the character associated to
it.
We shall consider different subgroups of $\mathfrak{S}_4$, each
one of them giving rise to a different equivariant model,
according to the following diagram (we use the following
shortenings: \texttt{GMM }for the general Markov model,
\texttt{K81} for the Kimura 3 parameter model, \texttt{K80} for
the Kimura 2 parameter model, \texttt{CS05} for the strand
symmetric model and \texttt{JC69} for the Jukes-Cantor model):
\begin{eqnarray*}
\tiny
 \xymatrix{\{\mathrm{id}\}  \ar[d] \ar[ddr] \\
{ \langle (\a\c)(\g\t),(\a\g)(\c\t)}\rangle \ar[d] \\
\langle (\a\c\g\t),(\a\g)\rangle \ar[d] & \langle (\a\t)(\c\g)\rangle \ar[ld] \\
\mathfrak{S}_4 &} & \tiny
\xymatrix{\mathtt{GMM}   \\
 \mathtt{K81} \ar[u]  \\
\mathtt{K80} \ar[u] & \mathtt{CS05} \ar[luu] \\
\mathtt{JC69} \ar[u] \ar[ru] &}
\end{eqnarray*}
Our aim here is to describe in a unified fashion the edge
invariants associated to these models for the case of a quartet
tree topology $T$, with leaves  $v_1,v_2,v_3,v_4$.
Write $e=L_1\mid L_2$ for the edge split corresponding to $e$, so
that $L_1=\{v_1,v_2\}$ and $L_2=\{v_3,v_4\}$.
\begin{figure}[here]
\begin{center}
\psfrag{v1}{$v_1$}\psfrag{v2}{$v_2$}\psfrag{v3}{$v_3$}
\psfrag{v4}{$v_4$} \psfrag{e}{$e$}
\includegraphics[scale=0.5]{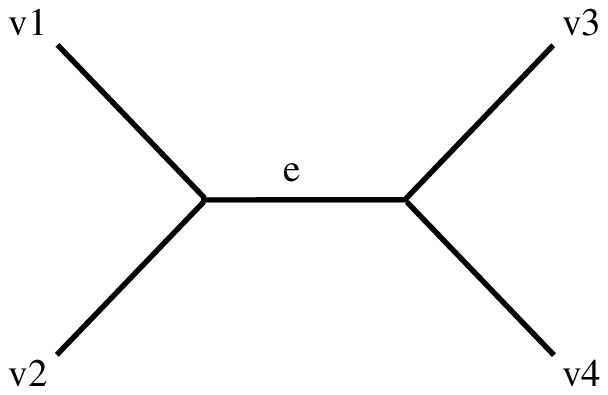}
\end{center}
\end{figure}

\begin{rem}\label{abelian}
When the subgroup $G\subset \mathfrak{S}_4$ is abelian, the usual
product of complex numbers induces on $\Omega_G$ a group
structure. Then, if $\{u^t_1,\ldots, u^t_{m_t}\}$ is a basis for
$W[\omega_t]$, for every $\omega_t\in \Omega_G$, we have that
\begin{eqnarray*}
\{u^{i_1}_{j_1}\otimes \dots \otimes u^{i_l}_{j_l} \mid
\omega_{i_1}\ldots \omega_{i_l}=\omega_t \}
\end{eqnarray*}
is a $\CC$-basis for $(\otimes^l W)[\omega_t]$.
\end{rem}

%
%
\subsection{General Markov model} As a first example, consider the trivial subgroup $\{\id\}\subset \mathfrak{S}_4$. The corresponding equivariant model is the \emph{general Markov model}, which is the most general model in the Felsenstein hierarchy (see Ch.4 in \cite{ASCB2005}). Invariants for this model have been studied by Allman and Rhodes in \cite{Allman2003,Allman2004b}. In this case, there is only one irreducible representation $\omega:G\rightarrow \CC$ defined by mapping $(\id)$ to $1$. The character table is
\begin{center}
  \begin{tabular}{ l || c }
 $\Omega_{(1)}$ & id   \\ \hline
$\omega$ & 1   \\ \hline $\chi$     & 4
  \end{tabular}
\end{center}
It follows that $\chi=4 \omega$. Keeping the notation introduced
in \ref{notation}, we have $\m=(4)$ and $W= W[\omega]\cong
N_{\omega}\otimes \CC^4$.

Now, for the case of four leaves, we have  $\chi^2=16\omega$ and
$\m(2)=(16)$. Then, the ideal $I_e$ is defined by the condition
\begin{eqnarray*}
\brk (M)\leq (4)
\end{eqnarray*}
where $M\in \Hom_G((W\otimes W)[\omega],(W\otimes W)[\omega])\cong
\Hom_{\CC}(\CC^{16}, \CC^{16})$ is a matrix of indeterminates
whose columns and rows are indexed by  the set $\{X_1\otimes
X_2\}_{X_1,X_2\in B}$.
The ideal $I_e$ obtained by imposing the above rank condition is
generated by $\binom{16}{5}\binom{16}{5}$ polynomials of degree~5.

\subsection{Strand symmetric model} Take  $G=\langle (\a\t)(\c\g) \rangle$, which
is  isomorphic to $\mathbb{Z}/2\mathbb{Z}$.
The equivariant matrices for this group have the following
structure:
\begin{eqnarray*}
 \left ( \begin{array}{cccc}
 a & b & c & d \\
e & f  & g & h \\
 h & g & f & e \\
 d & c & b  & a
 \end{array} \right )
\end{eqnarray*}
The equivariant model associated to $G$ is the \emph{strand
symmetric model} introduced in \cite{CS}.
There are two irreducible characters $\omega_1,\omega_2$, and the
character table  is
\begin{center}
  \begin{tabular}{ l || c |  c   }
$\Omega_G$ & $\rm{id}$ & $(\a\t)(\c\g)$   \\ \hline
$\omega_1$ & 1 & 1  \\
$\omega_2$ & 1 & -1 \\ \hline $\chi$     & 4 & 0
  \end{tabular}
\end{center}
Notice that since $G$ is abelian, all the irreducible
representations have dimension one. It follows that $\chi=2
\omega_1+2\omega_2$. Thus,  $\m=(2,2)$ and we have a decomposition
(\cite[Corollary 2.14]{FultHarr})
\[W= W[\omega_1]\oplus  W[\omega_2],\] where $W[\omega_{1}]\cong N_{\omega_1}\otimes \CC^2$ and  $W[\omega_{2}]\cong N_{\omega_2}\otimes \CC^2$.
Indeed, if we write
\begin{eqnarray*}
\u_1=\a+\t \qquad \u_2=\c+\g \qquad  \v_1=\a-\t \qquad \v_2=\c-\g,
\end{eqnarray*}
we have
\begin{eqnarray*}
W[\omega_{1}] =\langle \u_1,\u_2 \rangle_{\CC} \qquad
W[\omega_{2}]=\langle \v_1, \v_2 \rangle_{\CC}.
\end{eqnarray*}

Now, we focus on the case of the tree with four leaves. We have
$\chi^2=8\omega_1+8\omega_2$, so $\m(2)=(8,8)$.  Moreover, using
that $G$ is abelian (see Remark \ref{abelian})
\begin{eqnarray*}
W\otimes W[\omega_1] & = & \langle \u_1\otimes \u_1, \u_1\otimes
\u_2, \u_2\otimes \u_1, \u_2\otimes \u_2, \v_1\otimes \v_1,
\v_1\otimes \v_2, \v_2\otimes \v_1, \v_2\otimes \v_2
\rangle \\
W\otimes W[\omega_2]& = &\langle \u_1\otimes \v_1, \u_1\otimes
\v_2, \u_2\otimes \v_1, \u_2\otimes \v_2, \v_1\otimes \u_1,
\v_1\otimes \u_2, \v_2\otimes \u_1, \v_2\otimes \u_2
 \rangle
\end{eqnarray*}
Then, the ideal $I_e$ is defined by the conditions
\begin{eqnarray*}
\brk \left (\begin{array}{cccc}
M_1 & 0  \\
0 & M_2
  \end{array}
\right ) \leq (2,2)
\end{eqnarray*}
where
\begin{eqnarray*}
M_{1} = \left (\begin{array}{cccccccc}
q_{\tt u_1u_1u_1u_1} & q_{\tt u_1u_1u_1u_2} & q_{\tt u_1u_1u_2u_1} & q_{\tt u_1u_1u_2u_2} & q_{\tt u_1u_1v_1v_1} & q_{\tt u_1u_1v_1v_2} & q_{\tt u_1u_1v_2v_1} & q_{\tt u_1u_1v_2v_2 }\\
q_{\tt u_1u_2u_1u_1} & q_{\tt u_1u_2u_1u_2} & q_{\tt u_1u_2u_2u_1} & q_{\tt u_1u_2u_2u_2} & q_{\tt u_1u_2v_1v_1} & q_{\tt u_1u_2v_1v_2} & q_{\tt u_1u_2v_2v_1} & q_{\tt u_1u_2v_2v_2 }\\
q_{\tt u_2u_1u_1u_1} & q_{\tt u_2u_1u_1u_2} & q_{\tt u_2u_1u_2u_1} & q_{\tt u_2u_1u_2u_2} & q_{\tt u_2u_1v_1v_1} & q_{\tt u_2u_1v_1v_2} & q_{\tt u_2u_1v_2v_1} & q_{\tt u_2u_1v_2v_2 }\\
q_{\tt u_2u_2u_1u_1} & q_{\tt u_2u_2u_1u_2} & q_{\tt u_2u_2u_2u_1} & q_{\tt u_2u_2u_2u_2} & q_{\tt u_2u_2v_1v_1} & q_{\tt u_2u_2v_1v_2} & q_{\tt u_2u_2v_2v_1} & q_{\tt u_2u_2v_2v_2 }\\
q_{\tt v_1v_1u_1u_1} & q_{\tt v_1v_1u_1u_2} & q_{\tt v_1v_1u_2u_1} & q_{\tt v_1v_1u_2u_2} & q_{\tt v_1v_1v_1v_1} & q_{\tt v_1v_1v_1v_2} & q_{\tt v_1v_1v_2v_1} & q_{\tt v_1v_1v_2v_2 }\\
q_{\tt v_1v_2u_1u_1} & q_{\tt v_1v_2u_1u_2} & q_{\tt v_1v_2u_2u_1} & q_{\tt v_1v_2u_2u_2} & q_{\tt v_1v_2v_1v_1} & q_{\tt v_1v_2v_1v_2} & q_{\tt v_1v_2v_2v_1} & q_{\tt v_1v_2v_2v_2 }\\
q_{\tt v_2v_1u_1u_1} & q_{\tt v_2v_1u_1u_2} & q_{\tt v_2v_1u_2u_1} & q_{\tt v_2v_1u_2u_2} & q_{\tt v_2v_1v_1v_1} & q_{\tt v_2v_1v_1v_2} & q_{\tt v_2v_1v_2v_1} & q_{\tt v_2v_1v_2v_2 }\\
q_{\tt v_2v_2u_1u_1} & q_{\tt v_2v_2u_1u_2} & q_{\tt v_2v_2u_2u_1}
& q_{\tt v_2v_2u_2u_2} & q_{\tt v_2v_2v_1v_1} & q_{\tt
v_2v_2v_1v_2} & q_{\tt v_2v_2v_2v_1} & q_{\tt v_2v_2v_2v_2 }
  \end{array}
\right )\\
M_{2} =\left (\begin{array}{cccccccc}
q_{\tt u_1v_1u_1v_1} & q_{\tt u_1v_1u_1v_2} & q_{\tt u_1v_1u_2v_1} & q_{\tt u_1v_1u_2v_2} & q_{\tt u_1v_1v_1u_1} & q_{\tt u_1v_1v_1u_2} & q_{\tt u_1v_1v_2u_1} & q_{\tt u_1v_1v_2u_2 }\\
q_{\tt u_1v_2u_1v_1} & q_{\tt u_1v_2u_1v_2} & q_{\tt u_1v_2u_2v_1} & q_{\tt u_1v_2u_2v_2} & q_{\tt u_1v_2v_1u_1} & q_{\tt u_1v_2v_1u_2} & q_{\tt u_1v_2v_2u_1} & q_{\tt u_1v_2v_2u_2 }\\
q_{\tt u_2v_1u_1v_1} & q_{\tt u_2v_1u_1v_2} & q_{\tt u_2v_1u_2v_1} & q_{\tt u_2v_1u_2v_2} & q_{\tt u_2v_1v_1u_1} & q_{\tt u_2v_1v_1u_2} & q_{\tt u_2v_1v_2u_1} & q_{\tt u_2v_1v_2u_2 }\\
q_{\tt u_2v_2u_1v_1} & q_{\tt u_2v_2u_1v_2} & q_{\tt u_2v_2u_2v_1} & q_{\tt u_2v_2u_2v_2} & q_{\tt u_2v_2v_1u_1} & q_{\tt u_2v_2v_1u_2} & q_{\tt u_2v_2v_2u_1} & q_{\tt u_2v_2v_2u_2 }\\
q_{\tt v_1u_1u_1v_1} & q_{\tt v_1u_1u_1v_2} & q_{\tt v_1u_1u_2v_1} & q_{\tt v_1u_1u_2v_2} & q_{\tt v_1u_1v_1u_1} & q_{\tt v_1u_1v_1u_2} & q_{\tt v_1u_1v_2u_1} & q_{\tt v_1u_1v_2u_2 }\\
q_{\tt v_1u_2u_1v_1} & q_{\tt v_1u_2u_1v_2} & q_{\tt v_1u_2u_2v_1} & q_{\tt v_1u_2u_2v_2} & q_{\tt v_1u_2v_1u_1} & q_{\tt v_1u_2v_1u_2} & q_{\tt v_1u_2v_2u_1} & q_{\tt v_1u_2v_2u_2 }\\
q_{\tt v_2u_1u_1v_1} & q_{\tt v_2u_1u_1v_2} & q_{\tt v_2u_1u_2v_1} & q_{\tt v_2u_1u_2v_2} & q_{\tt v_2u_1v_1u_1} & q_{\tt v_2u_1v_1u_2} & q_{\tt v_2u_1v_2u_1} & q_{\tt v_2u_1v_2u_2 }\\
q_{\tt v_2u_2u_1v_1} & q_{\tt v_2u_2u_1v_2} & q_{\tt v_2u_2u_2v_1}
& q_{\tt v_2u_2u_2v_2} & q_{\tt v_2u_2v_1u_1} & q_{\tt
v_2u_2v_1u_2} & q_{\tt v_2u_2v_2u_1} & q_{\tt v_2u_2v_2u_2 }
\end{array}
\right )
\end{eqnarray*}
and $q_{xyzt}$ are the coordinates in the basis $x\otimes y\otimes
z\otimes t$.
We see that $I_e$ is generated by $\binom{8}{3}\binom{8}{3}+\binom{8}{3}\binom{8}{3}=6272$ polynomials of degree 3.

\subsection{Kimura 3-parameter model} Take $G=\langle (\a\c)(\g\t),(\a\g)(\c\t)\rangle$, which is also isomorphic to
$\mathbb{Z}/2\mathbb{Z} \times \mathbb{Z}/2\mathbb{Z}$.
The equivariant matrices for this group have the following
structure:
\begin{eqnarray*}
 \left ( \begin{array}{cccc}
 a & b & c & d \\
b & a  & d & c \\
 c & d & a & b \\
 d & c & b  & a
 \end{array} \right )
\end{eqnarray*}
In this case, the equivariant model is the \emph{Kimura
3-parameter model} introduced in \cite{Kimura1981}.
We write $\omega_{\a},\omega_\c,\omega_\g,\omega_\t$
for the irreducible characters of $G$.
The corresponding table is
\begin{center}
  \begin{tabular}{ l || c |  c |  c |  c  }
 $\Omega_G$ & $\textrm{id}$ & $(\a\c)(\g\t)$ & $(\a\g)(\c\t)$ & $(\a\t)(\c\g)$  \\ \hline
$\omega_\a$ & 1 & 1 & 1 & 1  \\
$\omega_\c$ & 1 & -1 & 1 & -1  \\
$\omega_\g$ & 1 & 1 & -1 & -1  \\
$\omega_\t$ & 1 & -1 & -1 & 1  \\ \hline $\chi$     & 4 & 0 & 0 &
0
  \end{tabular}
\end{center}
It follows that $\chi=\omega_\a+\omega_\c+\omega_\g+\omega_\t$ and
so, $\m=(1,1,1,1)$
\[W= W[\omega_\a]\oplus W[\omega_\c]\oplus W[\omega_\g]\oplus W[\omega_\t],\]where
\begin{eqnarray*}
W[\omega_{\a}]\cong N_{\omega_\a} \qquad W[\omega_{\c}]\cong
N_{\omega_\c}  \qquad W[\omega_{\g}]\cong N_{\omega_\g} \qquad
W[\omega_{\t}]\cong N_{\omega_\t}.
\end{eqnarray*}
In fact, if we write \begin{eqnarray}\label{q's}
\tiny & \overline{\a} & =\a+\c+\g+\t \qquad   \overline{\c}=\a+\c-\g-\t \\
\nonumber & \overline{\g} & =\a-\c+\g-\t \qquad
\overline{\t}=\a-\c-\g+\t
\end{eqnarray}
we have
\begin{eqnarray*}
W[\omega_{\a}]=\langle \overline{\a} \rangle \qquad
W[\omega_{\c}]=\langle \overline{\c} \rangle \qquad
W[\omega_{\g}]=\langle \overline{\g} \rangle \qquad
W[\omega_{\t}]=\langle \overline{\t} \rangle
\end{eqnarray*}
We remark that the basis
$\{\overline{\a},\overline{\c},\overline{\g},\overline{\t}\}$ is
the image of $\{\a,\c,\g,\t\}$ by the Fourier transform described
in \cite{CFS2} or \cite{CGS}.

Since  $\chi^2=4\omega_\a+4\omega_\c+4\omega_\g+4\omega_\t$, we
have  $\m(2)=(4,4,4,4)$. In virtue of Remark \ref{abelian},
\begin{eqnarray*}
W\otimes W[\omega_{\a}]=\langle  \overline{\a}\otimes
\overline{\a}, \overline{\c}\otimes \overline{\c},
\overline{\g}\otimes \overline{\g}, \overline{\t}\otimes
\overline{\t}
\rangle \\
W\otimes W[\omega_{\c}]=\langle  \overline{\a}\otimes
\overline{\c}, \overline{\c}\otimes \overline{\a},
\overline{\g}\otimes \overline{\t}, \overline{\t}\otimes
\overline{\g}
\rangle \\
W\otimes W[\omega_{\g}]=\langle  \overline{\a}\otimes
\overline{\g}, \overline{\c}\otimes \overline{\t},
\overline{\g}\otimes \overline{\a}, \overline{\t}\otimes
\overline{\c}
\rangle \\
W\otimes W[\omega_{\t}]=\langle  \overline{\a}\otimes
\overline{\t}, \overline{\c}\otimes \overline{\g},
\overline{\g}\otimes \overline{\c}, \overline{\t}\otimes
\overline{\a} \rangle
\end{eqnarray*}
Then, $I_e$ is given by the conditions
\begin{eqnarray}\label{rank3}
\brk \left (\begin{array}{cccc}
M_\a & 0 & 0 & 0  \\
0 & M_\c & 0 & 0  \\
0 & 0 & M_\g & 0  \\
0 & 0 & 0 & M_\t
  \end{array}
\right ) \leq (1,1,1,1)
\end{eqnarray}
where $M_{Z}\in M_{4,4}$ for all $Z\in B$, that is,
\begin{eqnarray*}
\tiny M_\a=\left (\begin{array}{cccc}
q_{\a\a\a\a} & q_{\a\a\c\c} & q_{\a\a\g\g} & q_{\a\a\t\t}  \\
q_{\c\c\a\a} & q_{\c\c\c\c} & q_{\c\c\g\g} & q_{\c\c\t\t}  \\
q_{\g\g\a\a} & q_{\g\g\c\c} & q_{\g\g\g\g} & q_{\g\g\t\t}  \\
q_{\t\t\a\a} & q_{\t\t\c\c} & q_{\t\t\g\g} & q_{\t\t\t\t}
  \end{array}
\right ) \quad \tiny M_\c=\left (\begin{array}{cccc}
q_{\a\c\a\c} & q_{\a\a\c\a} & q_{\a\a\g\t} & q_{\a\a\t\g}  \\
q_{\c\a\a\c} & q_{\c\a\c\a}& q_{\c\a\g\t} & q_{\c\a\t\g}  \\
q_{\g\t\a\c} & q_{\g\t\c\a} & q_{\g\t\g\t} & q_{\g\t\t\g}  \\
q_{\t\g\a\c} & q_{\t\g\c\a} & q_{\t\g\g\t} & q_{\t\g\t\g}
  \end{array}
\right ) \\
\tiny M_\g=\left (\begin{array}{cccc}
q_{\a\g\a\g} & q_{\a\g\c\t} & q_{\a\g\g\a} & q_{\a\g\t\c}  \\
q_{\c\t\a\g} & q_{\c\t\c\t} & q_{\c\t\g\a} & q_{\c\t\t\c}  \\
q_{\g\a\a\g} & q_{\g\a\c\t} & q_{\g\a\g\a} & q_{\g\a\t\c}  \\
q_{\t\c\a\g} & q_{\t\c\c\t} & q_{\t\c\g\a} & q_{\t\c\t\c}
  \end{array}
\right ) \quad \tiny M_\t=\left (\begin{array}{cccc}
q_{\a\t\a\t} & q_{\a\t\c\g} & q_{\a\t\g\c} & q_{\a\t\t\a}  \\
q_{\c\g\a\t} & q_{\c\g\c\g}& q_{\c\g\g\c} & q_{\c\g\t\a}  \\
q_{\g\c\a\t} & q_{\g\c\c\g} & q_{\g\c\g\c} & q_{\g\c\t\a}  \\
q_{\t\a\a\t} & q_{\t\a\c\g} & q_{\t\a\g\c} & q_{\t\a\t\a}
  \end{array}
\right )
\end{eqnarray*}
where $q_{\tt X_1X_2X_3X_4}$ are the coordinates in the basis
$\{\tt \overline{X}_1\otimes \overline{X}_2\otimes
\overline{X}_3\otimes \overline{X}_4\}_{X_i\in B}$.
The ideal $I_e$ obtained by imposing the rank conditions of
(\ref{rank3}) is generated by $\binom{4}{2}\binom{4}{2}+\binom{4}{2}\binom{4}{2}+\binom{4}{2}\binom{4}{2}+\binom{4}{2}\binom{4}{2}=144$ quadrics. However, at any point of $V(I_e)$ the variety is locally defined by 36 quadrics (see \cite[Example
4.9]{CFS2}).
\subsection{Kimura 2-parameter model} Take $G=\langle (\a\c\g\t),(\a\g) \rangle$, which is isomorphic to the dihedral group.
The equivariant matrices for this group have the following
structure:
\begin{eqnarray*}
 \left ( \begin{array}{cccc}
 a & b & c & b \\
b & a  & b & c \\

 c & b & a & b \\
 b & c & b  & a
 \end{array} \right )
\end{eqnarray*}
The equivariant model is the \emph{Kimura 2-parameter model}
introduced in \cite{Kimura1980}.
There are 5 irreducible characters $\omega_1$, $\omega_2$,
$\omega_3$, $\omega_4$, $\omega$ and the corresponding table is
\begin{center}
  \begin{tabular}{ l || c |  c |   c |   c |  c }
 $\Omega_{G}$ & $\rm{id}$ & $(\a\c\g\t)$ & $(\a\g)$ & $(\a\g)(\c\t)$ & $(\a\t\g\c)$ \\ \hline
$\omega_1$ & 1 & 1 & 1 & 1 & 1  \\
$\omega_2$ & 1 & 1 & -1 & 1 & 1  \\
$\omega_3$ & 1 & -1 & 1  & 1 & -1  \\
$\omega_4$ & 1 & -1 & -1  & 1 & -1 \\
$\omega$ & 2 & 0 & 0  & -2 & 0 \\ \hline $\chi$     & 4 & 0 & 2 &
0 & 0
  \end{tabular}
\end{center}
Notice that $G$ is not abelian and that the irreducible
representation $\omega$ is 2-dimensional. It follows that
$\chi=\omega_1+\omega_3+\omega$ and so, $\m=(1,0,1,0,1)$ and
\[W= W[\omega_1]\oplus W[\omega_3]\oplus W[\omega],\]where
\begin{eqnarray*}
W[\omega_{1}]\cong N_{\omega_1} \qquad W[\omega_{3}]\cong
N_{\omega_3}  \qquad W[\omega]\cong N_{\omega}.
\end{eqnarray*}
In fact, with the notation of (\ref{q's}) we have
\begin{eqnarray*}
W[\omega_1]=\langle \overline{\a} \rangle \qquad
W[\omega_3]=\langle \overline{\g} \rangle \qquad W[\omega]=\langle
\overline{\c}, \overline{\t}  \rangle\qquad
\end{eqnarray*}

Now, we consider the case of four leaves.
We have $\chi^2=3\omega_1+\omega_2+ 3\omega_3+\omega_4+4\omega$,
so $\m(2)=(3,1,3,1,4)$.
If $\psi\in \LL(T)^{G}$, then
\begin{eqnarray*}
\tiny \sf{e}{\psi} =  \left (\begin{array}{ccccc}
S_1 & 0 & 0 & 0 & 0  \\
0 & S_2 & 0 & 0  & 0  \\
0 & 0 & S_3   & 0 & 0  \\
0 & 0 & 0   & S_4 & 0 \\
0 & 0 & 0  &  0 & S
\end{array}
\right )\in M_{\m(2),\m(2)}
\end{eqnarray*}
where
\begin{eqnarray*}
S_1\in M_{3,3}\qquad S_2\in M_{1,1}\qquad S_3\in M_{3,3}\qquad
S_3\in M_{1,1}\qquad S\in M_{4,4}.
\end{eqnarray*}
Then, the ideal $I_e$ is given by the condition
\begin{eqnarray}\label{rank3}
\brk \sf{L_1\mid L_2}{\psi} \leq (1,0,1,0,1).
\end{eqnarray}
By imposing these rank conditions to the matrix
$\sf{L_1,L_2}{\psi}$ we obtain $\binom{3}{2}\binom{3}{2}+\binom{1}{1}\binom{1}{1}+
\binom{3}{2}\binom{3}{2}+\binom{1}{1}\binom{1}{1}+\binom{4}{2}\binom{4}{2}=9+1+9+1+36=56$  invariants: 54 of them are quadrics and 2 of them
are linear invariants.

\subsection{Jukes-Cantor model} Finally, we take the whole group of permutations $\mathfrak{S}_4$.
The equivariant matrices for this group have the following
structure:
\begin{eqnarray*}
 \left ( \begin{array}{cccc}
 a & b & b & b \\
b & a  & b & b \\
 b & b & a & b \\
 b & b & b  & a
 \end{array} \right )
\end{eqnarray*}
 The equivariant model associated to it is the \emph{Jukes-Cantor model} introduced in \cite{JC69}. The group $\mathfrak{S}_4$ has five irreducible characters $\{\omega_i\}_{i=0,\ldots,4}$ (see \S 2.3 of \cite{FultHarr}) and the following character table:

\begin{center}  \begin{tabular}{ l || c |  c |  c |  c | c  } $\Omega_{\mathfrak{S}_4}$ & $\mathrm{id}$ & $(\a\c)$ & $(\a\c\g)$ & $(\a\c\g\t)$ & $(\a\c)(\g\t)$ \\ \hline
$\omega_0$ & 1 & 1 & 1 & 1 & 1 \\
$\omega_1$ & 1 & -1 & 1 & -1 & 1 \\
$\omega_2$ & 2 & 0 & -1 & 0 & 2 \\
$\omega_3$ & 3 & 1 & 0 & -1 & -1 \\
$\omega_4$ & 3 & -1 & 0 & 1 &- 1 \\
\hline $\chi$     & 4 & 2 & 1 & 0 & 0
  \end{tabular}
\end{center}

\vspace{3mm}
It follows that \[\chi=\omega_0+\omega_3,\] that is, $\chi$ is the
sum of the trivial and the \emph{standard} representations. We
have $\m=(1,0,0,1,0)$.  Thus, there is a decomposition
\[W= W[\omega_0]\oplus W[\omega_3],\]
where
\begin{eqnarray*}
W[\omega_0]\cong  N_{\omega_0}\otimes \CC^{m_0}\cong  N_{\omega_0} & \quad & \dim W[\omega_0]=1\\
W[\omega_3]\cong  N_{\omega_3}\otimes \CC^{m_3}\cong  N_{\omega_3}
& \quad & \dim W[\omega_3]=3.
\end{eqnarray*}
In fact, with the notation of (\ref{q's}), we have
\begin{eqnarray*} W[\omega_0] = \langle \overline{\a}\rangle
\qquad W[\omega_3] = \langle
\overline{\c},\overline{\g},\overline{\t}\rangle.
 \end{eqnarray*}

The ideal $I_e$ is generated by the $(m_j+1)$-minors of the $j$-th
box of $\sf{e}{\psi}$ with $j=1,2,\ldots,5$.
On the other hand, it is straightforward to see that  $\chi^2=2\omega_0+\omega_2+3\omega_3+\omega_4$, so  $\m(2)=(2,0,1,3,1)$  and we have \begin{eqnarray*}(W\otimes W)[\omega_0]& = & \langle q_{\a\a}, q_{\c\c}+q_{\g\g}+q_{\t\t}\rangle \\(W\otimes W)[\omega_2]& = &\langle q_{\c\c}-q_{\g\g}, q_{\c\c}-q_{\t\t}\rangle \\(W\otimes W)[\omega_3]& = &\langle q_{\a\c},q_{\a\g}, q_{\a\t},q_{\c\a}, q_{\g\a}, q_{\t\a}, 
q_{\c\t}+q_{\t\c}, q_{\c\g}+q_{\g\c},q_{\g\t}+q_{\t\g}\rangle
\\
(W\otimes W)[\omega_4]& = &\langle q_{\c\t}-q_{\t\c},
q_{\c\g}-q_{\g\c}, q_{\g\t}-q_{\t\g}\rangle
\end{eqnarray*}
and $q_{\texttt{XY}}=q_{\texttt{X}}\otimes q_{\texttt{Y}}$, for
any $\texttt{X},\texttt{Y}\in B$.
Now, if $\psi\in \LL(T)^{\mathfrak{S}_4}$ we have
\begin{eqnarray*}
\tiny \sf{e}{\psi} =  \left (\begin{array}{cccc}
S_0 & 0 & 0 & 0  \\
0 & S_2 & 0 & 0  \\
0 & 0 & S_3 & 0  \\
0 & 0 & 0 & S_4
  \end{array}
\right )\in M_{\m(2),\m(2)}
\end{eqnarray*}
where
\begin{eqnarray*}
S_0\in M_{2,2}\qquad S_2\in M_{1,1} \qquad S_3\in M_{3,3} \qquad
S_4\in M_{1,1}.
\end{eqnarray*}
For instance, we have
\begin{eqnarray*}
S_{0}=\left ( \begin{array}{c|c} q_{\a\a\a\a} &
q_{\a\a\c\c}+q_{\a\a\g\g}+q_{\a\a\t\t} \\ \hline
q_{\c\c\a\a}+q_{\g\g\a\a}+q_{\t\t\a\a} & { \tiny  \begin{array}{c}q_{\c\c\c\c}+q_{\g\g\c\c}+q_{\t\t\c\c}+\\
q_{\c\c\g\g}+q_{\g\g\g\g}+q_{\t\t\g\g}+\\
q_{\c\c\t\t}+q_{\g\g\t\t}+q_{\t\t\t\t} \end{array}}
\end{array}\right )
\end{eqnarray*}
while
\begin{eqnarray*}\tiny
S_{2}=(q_{\c\c\c\c}-q_{\c\c\g\g}-q_{\g\g\c\c}+q_{\g\g\g\g}).
\end{eqnarray*}

Now, given $\psi\in \LL(T)^{\mathfrak{S}_4}$, we have $\psi\in
V(T)$ if and only if
\begin{eqnarray}\label{rank2}
\brk \sf{e}{\psi}\leq \m.
\end{eqnarray}
By imposing these rank conditions to the matrix
$\sf{e}{\psi}$ we obtain $\binom{2}{2}\binom{2}{2}+0+\binom{1}{1}\binom{1}{1}+ \binom{3}{2}\binom{3}{2}+\binom{1}{1}\binom{1}{1}=12$ phylogenetic invariants
$\{f_i\}_{i=1,\ldots,12}$:
\begin{itemize}
\item[1.] $f_1,\ldots,f_{10}$ have degree 2 and are obtained by the conditions ${\rk(S_0), \rk(S_3)=1}$
\item[2.]  $f_{11}, f_{12}$ have degree one and are obtained by the conditions $S_1,S_4=0$. These two invariants are equivalent to Lake's invariants (cf. \cite{Lake1987}).
\end{itemize}

\bibliography{invariants}
\bibliographystyle{alpha}
\end{document}